\def\draft{n}
\documentclass[11pt]{amsart}
\usepackage{fullpage,amssymb,epic,eepic,pb-diagram,lamsarrow,pb-lams,epsfig}

\theoremstyle{plain}

\newtheorem{theorem}{Theorem}
\newtheorem{proposition}{Proposition}[section]
\newtheorem{lemma}[proposition]{Lemma}
\newtheorem{corollary}[proposition]{Corollary}

\theoremstyle{definition}

\newtheorem{question}{Question}

\theoremstyle{remark}

\newtheorem{remark}[proposition]{Remark}

\def\printname#1{
	\if\draft y
		\smash{\makebox[0pt]{\hspace{-0.5in}
			\raisebox{8pt}{\tt\tiny #1}}}
	\fi
}

\newcommand{\psdraw}[2]
         {\begin{array}{c} \hspace{-1.3mm}
	\raisebox{-4pt}{\epsfig{figure=draws/#1.eps,width=#2}}
	\hspace{-1.9mm}\end{array}}

\newlength{\standardunitlength}
\catcode`\@=11
\long\def\@makecaption#1#2{%
     \vskip 10pt

\setbox\@tempboxa\hbox{
       \small\sf{\bfcaptionfont #1. }\ignorespaces #2}%
     \ifdim \wd\@tempboxa >\captionwidth {%
         \rightskip=\@captionmargin\leftskip=\@captionmargin
         \unhbox\@tempboxa\par}%
       \else
         \hbox to\hsize{\hfil\box\@tempboxa\hfil}%
     \fi}
\font\bfcaptionfont=cmssbx10 scaled \magstephalf
\newdimen\@captionmargin\@captionmargin=2\parindent
\newdimen\captionwidth\captionwidth=\hsize
\catcode`\@=12

\newcommand{\im}{\operatorname{Im}}

\def\lbl#1{\label{#1}\printname{#1}}

\def\eqdef{\overset{\text{def}}{=}}

\def\BZ{\mathbb Z}
\def\BQ{\mathbb Q}

\def\A{\mathcal A}

\def\G{\mathcal G}

\def\T{\mathcal T}
\def\M{\mathcal M}
\def\F{\mathcal F}

\def\T{\mathcal T}

\def\LL{\mathsf{L}}
\def\LLc{\mathsf{D}}
\def\Lie{\mathsf{D}'}
\def\sfY{\operatorname{\mathsf{Y}}}

\def\La{\Lambda}
\def\l{\lambda}
\def\Ga{\Gamma}
\def\S{\Sigma}
\def\si{\sigma}
\def\varsi{\varsigma}

\def\ihs{integral homology 3-sphere}

\def\fti{finite type
invariant}

\def\mcg{mapping class group}
\def\nequiv{$n$-equivalence}
\def\cok{\mathrm{cok}}
\def\otau{\overline{\tau}}
\def\torsion{\mathrm{torsion}}

\def\la{\langle}
\def\ra{\rangle}

\def\o1o{\underset{1}\ast}
\def\x1y{\underset{x1y}\ast}
\def\y1x{\underset{y1x}\ast}
\def\sti{\Sigma_{g,1}\times I}
\def\stig{\Sigma_{g}\times I}
\def\wb{W^{\mathfrak{b}}}

\def\FY#1{\mathcal F^Y_{#1}\mathcal M}

\def\GY#1{\mathcal G^Y_{#1}\mathcal M}

\def\s{$ \spadesuit $}

\def\we{\wedge}

\def\i{^{-1}}

\def\a{\alpha}

\def\bd{\partial}

\def\e{\epsilon}
\def\Ga{\Gamma}
\def\d{\delta}
\def\b{\beta}

\def\s{\sigma}

\def\sub{\subseteq}

\def\Sg{\Sigma_{g,1}}
\def\H{\mathcal H}
\def\Kk{K(F/F_{n},1)}
\def\eb{\mathfrak b}
\renewcommand{\ker}{\operatorname{Ker}}
\newcommand{\rk}{\operatorname{rank}}
\newcommand{\tor}{\operatorname{tors}}
\newcommand{\con}{\equiv}

\def\iso{\cong}
\def\pt{\mathbb Z[\![t_1 ,\dots,t_m ]\!]}
\def\sms{\smallsmile}       
\def\smf{\smallfrown}       

\begin{document}


\title[Trees, 3-manifolds, Massey products and Johnson's
homomorphism]{Tree-level invariants of three-manifolds,
Massey products and the Johnson homomorphism}

\author{Stavros Garoufalidis}
\address{School of Mathematics \\
         Georgia Institute of Technology \\
         Atlanta, GA 30332-0160, USA \\ 
         {\tt http://www.math.gatech} \newline {\tt .edu/$\sim$stavros } }
\email{stavros@math.gatech.edu}

\thanks{The  authors were partially supported by NSF grants
        DMS-98-00703 and DMS-96-26639 respectively, and by an Israel-US BSF
grant.\newline
        This and related preprints can also be obtained at
{\tt http://www.math.gatech.edu/$\sim$stavros } and
{\tt http://www.math.brandeis.edu/Faculty/jlevine/ }
\newline
1991 {\em Mathematics Classification.} Primary 57N10. Secondary 57M25.
\fti s, homology cylinders,
Johnson's homomorphism, deformation quantization, graph cohomology.
}

\author{Jerome Levine}
\address{Department of Mathematics\\
         Brandeis University\\
         Waltham, MA 02254-9110, USA. }
\email{levine@brandeis.edu}

\dedicatory{Dedicated to Dennis Sullivan on the occasion of his 60th
birthday}

\date{
July 25, 2003 \hspace{0.5cm} First edition: March 8,
1999.}

\begin{abstract}
We show that the tree-level part of a theory of finite type invariants of 
3-manifolds (based on surgery on objects called claspers, Y-graphs or 
clovers) is essentially given by classical algebraic topology
in terms of the Johnson homomorphism and Massey products, for arbitrary 
3-manifolds. A key role of our proof is played by the notion of a homology
cylinder, viewed as an enlargement of the mapping class group, and  an
apparently new Lie algebra of graphs colored by $H_1(\Sigma)$ of a closed 
surface $\Sigma$, closely related to deformation quantization on a surface
\cite{AMR1, AMR2,Ko3} as well as to a Lie algebra that encodes the symmetries 
of Massey products and the Johnson homomorphism. In addition, we give a 
realization theorem for Massey products and the Johnson homomorphism by 
homology cylinders.
\end{abstract}

\maketitle



\section{Introduction}
\lbl{sec.intro}
\subsection{A brief summary} In this paper we investigate relations 
between three different phenomena in low-dimensional topology:
\begin{itemize}
\item[(a)] Massey products on the first cohomology $H^1 (M)$ 
with integer coefficients of $3$-manifolds $M$.
\item[(b)] the Johnson homomorphism on the mapping class group of an 
orientable surface
\item[(c)] 
the Goussarov-Habiro theory of finite-type invariants of 3-manifolds.
\end{itemize}

A key point of the connection between (a) and (b) is the notion of a 
{\em homology cylinder}, i.e., a homology cobordism between 
an orientable surface and itself. This notion generalizes the mapping
class group of a surface (in that case the cobordism is a product).
We will construct an extension of the Johnson homomorphism to homology 
cylinders and use it to completely determine, in an 
explicit fashion, the possible Massey products at the first 
non-trivial level in a closed $3$-manifold (assuming that the first
homology $H_1$ is 
torsion-free)---see Theorem \ref{thm.1}, Corollary \ref{cor.1} and 
Theorem \ref{thm.2}. 

This generalizes the known relationship 
between the Johnson homomorphism and Massey products in the mapping 
torus of a diffeomorphism of a surface to the more general situation 
of homology cylinders---see Theorems \ref{thm.3},  \ref{th.alf} and 
Remark \ref{rem.hcc}.

For historical reasons, we should mention early work of Sullivan \cite{Su}
on a relation between (a) and (b), and, for an alternative
point of view,  work of Turaev \cite{Tu}.
 
With regards to the connection between (b) and (c), the main idea is
to consider Massey products as finite-type invariants of 
$3$-manifolds, and to interpret them by a graphical calculus on
{\em trees}---see
Theorem \ref{thm.7}---in much the same way that 
Vassiliev invariants of links have a graphical representation
and that Milnor's $\mu$-invariants are known to be 
exactly the Vassiliev invariants of (concordance classes of) 
string links which are represented by trees, see \cite{HM}.
A by-product of this investigation is a curious Lie algebra 
structure on a vector space of the graphs which describe  \fti s 
of homology cylinders---see Proposition \ref{thm.4} and 
Theorem \ref{thm.5}---that corresponds to the stacking of one homology
cylinder on top of another, and is closely related to
deformation quantization on a surface \cite{AMR1, AMR2,Ko3}.

\subsection{History}
\lbl{sub.his}
Years ago, Johnson introduced a homomorphism (the so-called Johnson
homomorphism) which he used to study the mapping class group, 
\cite{Jo1,Mo3}. Morita \cite{Mo} discovered a close relation between
the Johnson homomorphism and the simplest \fti\ of 3-manifolds, namely
the Casson invariant; this relation was subsequently generalized
by the authors \cite{GL1, GL2} to
all \fti s of \ihs s (i.e., 3-manifolds $M$ with $H_1(M,\BZ)=0$).
This generalization posed the question of
understanding  the Johnson homomorphism (crucial to the structure
of the mapping class group) from the point
of view of \fti s.  Unfortunately, this question is rather hard to answer
if we confine ourselves to invariants of \ihs s. This difficulty is overcome
by using a theory of \fti s based on the notion of surgery on $Y$-links, see 
\cite{Gu1,Gu2,Ha,Oh,GGP}. Using this theory we will  show that the Johnson 
homomorphism is contained in its tree-level part, and we conjecture
that an extension of the Johnson homomorphism to homology cylinders
(i.e., 3-manifolds with boundary that homologically look like the product
of a surface with $[0,1]$), which
we define below, gives the full tree-level part; thus answering
questions raised by Hain and Morita \cite{Hain,Mo3}.

En route to answering the above  question,
we were led to study this theory of invariants for
homology cylinders (studied also from a slightly different perspective
by Goussarov \cite{Gu1,Gu2} and Habiro \cite{Ha})
and discovered an apparently new Lie algebra of graphs colored by
$H_1(\Sigma)$
of a closed surface $\Sigma$, closely related to
deformation quantization on a surface \cite{AMR1, AMR2},
and to the curious  graded group
$\LLc(A)\eqdef \text{Ker}(A \otimes \LL(A) \to \LL(A))$,
(where $\LL(A)$ denotes the free Lie ring of a torsion-free abelian
group
$A$) studied independently by several authors with a variety of
motivations
\cite{Jo2,Mo1,Ih,Dr,Ko1,Ko2,O1,O2,HM}.

It turns out that Massey products of 3-manifolds naturally take
values in $\LLc(A)$,
and so does the Johnson homomorphism, which is also  closely related to
Massey products--- a fact well-known to Johnson \cite{Jo2}, and later
proved by Kitano \cite{Ki}.  However it is now known that the Johnson
homomorphism cannot
realize all elements of $\LLc(A)$, but we will see that one can achieve
this realizability by replacing
surface diffeomorphisms by homology cylinders.

The generalized Johnson homomorphism actually provides
universally-defined invariants of homology cylinders, lifting the only
partially-defined  Massey products. These are our explicit
  candidates for the full
tree-level part of the Goussarov-Habiro theory (for homology cylinders).
This phenomenon was already observed when
one replaces 3-manifolds by string-links up to homotopy, see Bar-Natan
\cite{B-N} or by string-links up to concordance, see Habegger-Masbaum
\cite{HM}.
On the other hand, Massey products apply to more general manifolds and they
should provide partially-defined finite-type invariants.

\section{Statement of the results}
\lbl{sec.res}

\subsection{Conventions}
\lbl{sub.conv}

$F$ will always stand for a free group and $H$ for a torsion-free
abelian
group. The {\em lower central series} of
a group $G$ is inductively defined by $G_1=G$ and $G_{n+1}=[G,G_{n}]$.
A group homomorphism $p: K \to G$ is called an $n$-{\em equivalence}
if
it induces an isomorphism $K/K_n \cong G/G_n$. All manifolds will be
oriented, and all maps between them
will preserve orientation, unless otherwise mentioned.
The boundary of an oriented manifold is oriented with the ``outward
normal
first'' convention.

\subsection{Massey products}
\lbl{sub.massey}

By {\em Massey products} of length $n \geq 2$ in $H^2(\pi)$ we mean a Massey
product
$\la a_1, \dots, a_n \ra \in H^2(\pi)$ for $a_i \in H^1(\pi)$,
\cite{Ma,FS},
which are defined assuming that the ones of length $n-1$ are defined
and
vanish. We have the following theorem on universal Massey products:

\begin{theorem}
\lbl{thm.1}
(i) Given a connected topological space $X$ and $2$-equivalence $p:
F \to \pi\eqdef\pi_1(X)$, then $X$ has
vanishing Massey products of length less than $n$ if and only
if $p$ is an \nequiv .
\newline
(ii) In that case, we have a short exact sequence\footnote{
Note that $\pi_n$ denotes the $n$th commutator subgroup of $\pi_1(X)$
and not the $n$th homotopy group of $X$.}
$$
H_2(X,\BZ) \to \LL_n(H_1(X,\BZ)) \to \pi_{n}/\pi_{n+1} \to 0,
$$
where the first map determines and is determined by
all length $n$ Massey
products (for a precise expression, see Corollary \ref{cor.fi}) and
the second is induced by the Lie bracket. \newline
(iii) In addition, we have that
$$\a_1\sms \la \a_2 ,\dots ,\a_{n+1}\ra =\la \a_1 ,\dots, a_n
\ra \sms\a_{n+1} \in H^2(\pi), $$
for any $\a_1 ,\dots ,\a_{n+1}\in H^1 (\pi )$.
\end{theorem}

\begin{remark}
For the dependence of the short exact sequence in the above theorem
on the
map $p$, see Remark \ref{rem.dep}.
If $X$ satisfies the hypothesis of Theorem \ref{thm.1}
we have dually, over $\BQ$:
$$
0 \to (\pi_{n}/\pi_{n+1})^\ast_{\BQ} \to \LL_n(H^1(X,\BQ)) \to H^2(X,
\BQ).
$$
Note that the first part of Theorem \ref{thm.1} appears in
\cite[Lemma 16]{O2}, and that the exact sequence above was first
suggested by Sullivan in \cite{Su} for $n=2$, and subsequently proven
by
Lambe in \cite{La} for $n=2$ using  different techniques
involving minimal models.
\end{remark}

\begin{corollary}
\lbl{cor.1}
Given an \nequiv\
$p :F \to \pi\eqdef\pi_1(M)$, where $M$ is a closed 3-manifold,
we have the exact sequence
$$
  H^\ast \to \LL_n(H) \to \pi_{n}/\pi_{n+1} \to 0$$
where $H = H_1(M,\BQ)$.
If $\mu_n(M,p) \in H \otimes \LL_n(H)$
denotes the first map (abbreviated by $\mu_n(M)$ if $p$ is clear),
then we have that $$ \mu_n(M,p) \in \LLc_n(H).$$
In particular, $\mu_n(M,p)=0$ if and only if $p: F\to \pi$ is an
$(n+1)$-equivalence.
\end{corollary}

Given an integer $n$ and a  torsion-free abelian group $H$,
it is natural to ask which elements
of $\LLc_n(H)$ are realized by 3-manifolds as above.  For this see
Theorem \ref{thm.2} below.

\subsection{The Johnson homomorphism}
\lbl{sub.johnson}

We now discuss the relation between Massey products and the Johnson
homomorphism.

Let $\Ga_{g,1}$ denote the  \mcg\ of a surface $\S_{g,1}$
of genus $g$ with one boundary component (i.e., the group of
  surface diffeomorphisms that pointwise preserve the
boundary),
and let $\Ga_{g,1}[n]$ denote its subgroup that consists of
surface diffeomorphisms that induce the identity on $\pi/\pi_{n+1}$,
where $\pi=\pi_1(\S_{g,1})$. In \cite{Jo1}, Johnson defined a
homomorphism
$$ \tau_n: \Ga_{g,1}[n] \to \LLc_{n+1}(H),$$
where $H=H_1(\S_{g,1},\BZ)$, which he further extended to the case of
a closed surface.
Johnson was well-aware of the relation between his
homomorphism and Massey products on mapping torii, i.e., on twisted
surface
bundles over a circle; see \cite[p. 171]{Jo2}, further
elucidated by Kitano \cite{Ki}. In the present note, we extend this
relation to Massey products that come from an arbitrary pair
$(\S,M)$ of an imbedding
$\iota: \S \hookrightarrow M$ of a closed surface (not necessarily
separating)
in a 3-manifold.
Fix a closed 3-manifold $M$ and an $(n+1)$-equivalence
$F\to \pi\eqdef\pi_1(M)$. Given a pair $(\S,M)$, and
$\phi \in \Ga[n]$, let $M_\phi$ denote the result of cutting $M$ along
$\S$, twisting by the element $\phi$ of its \mcg\ and gluing back.
In this case, there exists a canonical  cobordism $ N_\phi $ between $M$ and
$M_\phi$ such that the maps $\pi_1(M) \rightarrow \pi_1(N_\phi) \leftarrow
\pi_1(M_\phi)$ (induced by the inclusions $M,M_\phi \hookrightarrow
N_\phi$)
are $(n+1)$-equivalences; thus by Theorem \ref{thm.1}, $M_\phi$
has vanishing Massey products of length less than $n+1$. The ones of
length
$n+1$ on $M_{\phi}$ are determined in terms of those of $M$ and
the Johnson homomorphism as follows:

\begin{theorem}
\lbl{thm.3}
With the above assumptions, we have
$$ \mu_{n+1}(M_\phi)=\mu_{n+1}(M) + \iota_\ast \tau_n(\phi).$$
\end{theorem}

See also Remark \ref{rem.hcc}.

\subsection{Homology cylinders and realization}
\lbl{sub.homcyl}

It is well known \cite{Mo2,Hain} that the Johnson homomorphism
$\tau_n$
is not onto, in other words not every element of $\LLc_{n+1}(H)$ can
be realized by surface diffeomorphisms.  Generalizing
surface diffeomorphisms to a more general notion of homology
cylinders (defined below) allows
us to define an {\em ungraded}  version
of the Johnson homomorphism, which then
induces, on the associated graded level,  generalizations of the Johnson
homomorphisms. We will show that all of these are onto,
see Theorem \ref{th.alf}.  As an application of this result, we will show
that we can realize
every element in $\LLc_n(H)$ by 3-manifolds as in Corollary
\ref{cor.1}
and, in addition give a proof, free of spectral sequences, of
the isomorphism \eqref{eq.h3}, as mentioned above.

Let $\Sg$ denote the compact orientable surface of genus $g$  with
one boundary component. A {\em homology cylinder} over $\Sg$ is a
compact orientable $3$-manifold $M$ equipped with two
imbeddings $i^{-},i^{+}:\Sg\to\bd M$ so that $i^{+}$ is
orientation-preserving and $i^{-}$ is orientation-reversing and if
we denote $\S^{\pm}=\im i^{\pm}(\Sg )$, then $\bd M=\S^+\cup\S^-$ and
$\S^{+}\cap\S^{-}=\bd\S^{+}=\bd\S^{-}$. We also require that
$i^{\pm}$ be homology isomorphisms. We can multiply two homology
cylinders by identifying $\S^{-}$ in the first with $\S^{+}$ in
the second via the appropriate $i^{\pm}$. Thus $\H_{g,1}$, the set
of orientation-preserving diffeomorphism classes of homology
cylinders over $\Sg$ is a semi-group with an obvious identity.

There is a canonical homomorphism $\Ga_{g,1}\to\H_{g,1}$ that
sends $\phi$ to $(I\times\Sg, 0\times\text{id}, 1\times \phi)$.
Nielsen showed that the natural map $\Ga_{g,1}\to \mathrm{A}_0(F)$
is an isomorphism, where $F$ is
the free group on $2g$ generators $\{ x_i ,y_i\}$, identified with
the fundamental group of $\Sg$ (with base-point on $\bd\Sg$), and
$\mathrm{A}_0 (F)$ is the group of automorphisms of $F$ which fix the
element
$\omega_g =[x_1 ,y_1 ]\cdots [x_g ,y_g ]$, representing the boundary
of
$\Sg$. It is natural to ask whether there exists an analogous isomorphism
for the semigroup $\H_{g,1}$. Below, we construct for every $n$
a homomorphism $\si_n :\H_{g,1}\to \mathrm{A}_0
(F/F_n )$, where $\mathrm{A}_0 (F/F_n )$ is the group of
automorphisms $\phi$
of $F/F_n$ such that a lift of $\phi$ to an endomorphism $\bar \phi$
of $F$
fixes $\omega_g\mod F_{n+1}$. It is easy to see that this condition
is independent of the lift. For example $A_0 (F/F_2 )=\text{Sp}(g, \BZ )$.

Given $(M,i^{+},i^{-})\in\H_{g,1}$ consider
the homomorphisms $i^{\pm}_{*}:F\to\pi_{1} (M)$, where the
base-point is taken in $\bd\S^{+}=\bd\S^{-}$. In general, $i^{\pm}_{*}$
are not isomorphisms--- however,  since $i^{\pm}$ are
homology isomorphisms, it follows from Stallings \cite{St} that they
  induce isomorphisms
$i^{\pm}_{n}:F/F_{n}\to\pi_{1}(M)/\pi_{1}(M)_{n}$. We then
define $\si_{n}(M,i^{\pm})=(i^{-}_{n})^{-1} \circ i_{n}^{+}$. It is
easy to see that $\si_n(M,i^{\pm}) \in \mathrm{A}_0(F/F_n)$.

\begin{theorem}
\lbl{th.alf}
The map $\si_{n}: \H_{g,1}\to \mathrm{A}_0(F/F_n)$ is surjective.
\end{theorem}

\begin{remark}
\lbl{rem.hc}
We can convert $\H_{g,1}$ into a group $\H^c_{g,1}$ by considering
{\em homology cobordism classes } of homology cylinders. The inverse
of an
element is just the reflection in the $I$ coordinate. It is easy to
see
that the invariants $\si_n$ just depend on the homology bordism class
and so
define homomorphisms $\H^c_{g,1}\to\mathrm{A}_0 (F/F_n )$. The natural
homomorphism
$\Ga_{g,1}\to\H^c_{g,1}$ is seen to be injective by the existence of
the $\si_n$
and the fact that the homomorphism $\Ga_{g,1}\to\mathrm{A}_0 (F)$ is
an isomorphism.

In addition, we can combine the maps $\si_{n}$, for all $n$, to a single map
$\si^{\text{nil}}:\H_{g,1}\to \mathrm{A}_0(F^{\text{nil}})$,
where $F^{\text{nil}}$ is the nilpotent completion of $F$. Unlike
$\si$,
$\si^{\text{nil}}$ is not one-to-one, i.e.,
$\cap_{n}\ker\si_{n}\not=\{ 1\}$.
For example, if $P$ is any homology sphere, then the connected
sum $(I\times\S_{g,1} )\sharp P$ defines an element in the kernel.
Also
$\si^{\text{nil}}$ is not onto, even though each $\si_{n}$ is. To
identify the image of $\si^{\text{nil}}$
we have to consider the {\em algebraic closure }
$\bar F\sub F^{\text{nil}}$, see \cite{Le2}.  Using the arguments of
\cite{Le2}, we can show that any element of
$\im (\si^{\text{nil}})$ restricts to an automorphism of $\bar F$
and, by
  arguments similar to the proof of Theorem
\ref{th.alf}, it can be proved that $\im (\si^{\text{nil}})$ consists
precisely of those  $\phi\in A_0 (F^{\text{nil}})$ which restrict to an
automorphism of $\bar F$ and such that the element of $H_{2}(\bar F )$
associated to $\phi$ (see the
proof of Theorem \ref{th.alf}) is zero.
But since we do not know whether $H_{2} (\bar F
)=0$, this result does not seem very useful at this time.
\end{remark}

\begin{remark}\lbl{rem.mu}
The $\{\si_n\}$ can be described by numerical invariants if we consider the
coefficients of the Magnus expansion of $\si_n (M)(x_i ),\si_n (M) (y_i )$.
This is analogous to the definition of the $\mu$-invariants of a string
link. We can refer to these as {\em $\mu$-invariants of homology cylinders}.
\end{remark}

  It will be useful for us to consider the filtration  defined by the
maps $\si_n$, namely we define a decreasing {\em weight
filtration} on
$\H_{g,1}$ and on $\H^c_{g,1}$ by setting $\H_{g,1}[n]=\ker (\si_n)$.

\begin{proposition}
\lbl{prop.DA}
We have an exact sequence
$$ 1\to \LLc_n(H)\to \mathrm{A}_0(F/F_{n+1})\to \mathrm{A}_0(F/F_n)
\to 1$$
and a commutative diagram
$$
\begin{diagram}
\node{\Ga_{g,1}
[n]}\arrow[2]{e}\arrow{se,b}{\tau_n}\node[2]{\H_{g,1}[n]}
\arrow{sw,b}{\varsi_n}\\
\node[2]{\LLc_n(H)}
\end{diagram}
$$
where  the map $\varsi_n$, induced by $\s_n$, is onto. It follows that
$$ 0 \to \H^c_{g,1}[n+1] \to \H^c_{g,1}[n] \overset{\varsi_n}\to \LLc_n(H)
\to 0$$
is exact.
\end{proposition}
\noindent

\begin{remark}
A major problem in the study of the mapping class group is to
determine
the image of the Johnson homomorphism $\tau_n$, which largely
determines
the algebraic structure of the mapping class group since
$\cap_n \Ga_{g,1}[n]=1$. In contrast,  Theorem \ref{th.alf} largely
determines
the structure of $\H_{g,1}/\H_{g,1}[\infty]$,  but in this case
  $\H_{g,1}[\infty]=
\cap_n \H_{g,1}[n]$ is  not trivial---see Question \ref{q.inf} at the end
of the paper.
\end{remark}

\begin{remark}\lbl{rem.sl}
It is instructive to consider the analogy between, on the one hand, the
mapping class group,  homology cylinders and the invariant
$\s_n$ and the Johnson homomorphism, and, on the other hand, the pure braid
group,  string links and the Milnor $\mu$-invariants.
There is an injection of the pure braid group on $g$ strands into the mapping
class group $\Gamma_{g,1}$, first defined by Oda and studied in \cite{Le2},
which preserves the weight filtrations and induces a monomorphism of the
associated graded Lie algebras.
This can, in fact, be generalized to an injection of the semi-group
$\mathcal S_g$ of string
links on $g$ strands  into the semi-group $\H_{g,1}$ (and of the
string-link concordance group $ \mathcal S_g^c$ into $\H_{g,1}^c$), under
which $\s_n$ and
the $\mu$-invariants correspond. We will explain this in a future paper.

\end{remark}

\begin{theorem}
\lbl{thm.2}
Every element in $\LLc_n(H)$ is realized by an \nequiv\ $F\to
\pi_1(M^3)$, for some closed $3$-manifold $M$,
as in Corollary \ref{cor.1}. In addition, a $2$-equivalence $F \to H$
gives rise to a map $H_3(F/F_n)\to H \otimes \LL_n(H)$ inducing
the  isomorphism of Equation \eqref{eq.h3}.
\end{theorem}

We will give two different proofs of this theorem. One approach is to
apply results of Orr
\cite{O1} and Igusa-Orr \cite{IO} on $H_3 (F/F_n )$ and, in particular, the
isomorphism
\begin{equation}
\lbl{eq.h3}
\cok(H_3(F/F_{n+1}) \to H_3(F/F_n)) \cong \LLc_n(H).
\end{equation}
A very similar argument appears in \cite{CGO}.

Alternatively we will see that this realizability is a consequence of
Theorem \ref{th.alf}.
This approach has the advantage of being ``spectral-sequence-free'' and
also gives another proof of \eqref{eq.h3}.

\subsection{Homology cylinders and \fti s of 3-manifolds}
\lbl{sub.fti}

Goussarov and Habiro \cite{Gu1,Gu2,Ha} have studied
two rather dual notions: an $n$-equivalence relation among
3-manifolds,
and a theory of invariants of 3-manifolds with values in an abelian
group.
Since their work is recent and not yet fully written, we will,
for the benefit of the reader, give a short introduction using
terminology and notation from \cite{GGP} (to which we refer the reader
for detailed proofs).
Both notions are intimately related to that of surgery $M_{\Ga}$ along
a $Y$-link $\Ga$ in a 3-manifold $M$, i.e.,  surgery along an imbedded link
associated to an imbedding of an
appropriately oriented, framed graph with trivalent and univalent vertices
so that the
univalent
ones end in ``leaves'' (explained below). Two manifolds are $n$-{\em
equivalent} if one can pass from
one to the other by surgery on a $Y$-link  associated to a {\em
connected }  Y-graph of degree (i.e., number of
trivalent
vertices) at least $n$.  For example a theorem of Matveev \cite{M} says
that two closed manifolds are $1$-equivalent if and only if there is an
isomorphism between their first homology groups which preserves the
linking form on the torsion subgroups.
Similarly, a {\em finite type invariant} $\lambda$
ought to be the analog of a polynomial on the set of 3-manifolds,
in other words for some integer $n$ it satisfies a difference
equation
$$
\sum_{\Ga' \subseteq \Ga}(-1)^{|\Ga'|} \lambda(M_{\Ga})=0
$$
where $\Ga$ is a $Y$-link in $M$ of more than $n$
components and the sum is over all $Y$-sublinks $\Ga'$ of $\Ga$.
In view of the above definition, it is natural to consider the free abelian
group $\M$  generated by homeomorphism classes of closed oriented
$3$-manifolds, and to define a decreasing filtration $\FY {}$
on $\M$ in such a way that $\lambda$ is an invariant of type $n$ if and
only if it vanishes on $\FY {n+1}$. Thus the question of how many
invariants
of degree $n$ there are translates into a question about the size of
the graded quotients $\GY {n}\eqdef\FY {n}/\FY {n+1}$. One
traditionally
approaches this problem by giving independently an upper bound and a
lower bound, which hopefully match.
In this theory, an upper bound has been obtained in terms of an
abelian group
of decorated graphs as follows. One observes
first that surgery along $Y$-links preserves the homology and linking
form of 3-manifolds,  as well as the boundary. Define an equivalence
relation on compact $3$-manifolds: $M\sim N$ if there exists an isomorphism
$\rho :H_1 (M)\to H_1 (N)$ inducing an isometry of the linking forms, and a
homeomorphism $\bd M\to\bd N$ consistent with $\rho$. Thus if we let
$\M(M)$ denote the subgroup of
$\M$ generated by equivalent 3-manifolds, we have a direct sum
decomposition $\M=\oplus_{\sim}\ \M(M)$
(and also, $\FY {}=\oplus_{\sim}\ \FY {}(M)$), where the sum is over a
choice of one manifold $M$ from each equivalence class. In fact for closed
$3$-manifolds Matveev's theorem tells us that $ \GY {0}(M)=\BZ$. After we fix
a 3-manifold $M$, and an oriented link
$\mathfrak{b}$
in $M$ that represents a basis of $H_1(M,\BZ)/\torsion$, together with
a framing of $\mathfrak{b}$ (i.e., a choice of a trivialization of the normal
bundle of each component of $\mathfrak{b}$),
it turns out that there is a map\footnote{a more precise notation,
which
we will not use, would be $W^{M,\mathfrak{b}}_n$.}
\begin{equation}
\lbl{eq.graded}
\wb_n: \A_n(M) \to \GY {n}(M),
\end{equation}
which is onto over $\BQ$ (actually, onto over $\BZ [1/(2 |\torsion|)]$),
where $\A(M)$ is the group generated by graphs with univalent and
trivalent vertices, with a cyclic order along each trivalent
vertex, decorated by an element of $H_1(M,\BZ)$  on each univalent
vertex,
modulo some relations, see \cite{GGP}.
Here $\A_n(M)$ is the subgroup generated by
graphs of degree $n$, i.e., with $n$ trivalent vertices; thus we have
$\A(M)=\oplus_n\A_n (M)=\A^t(M) \oplus \A^l(M)$,
where $\A^t(M)$ (resp. $\A^l(M)$) is the subgroup of $\A(M)$
generated by trees (resp. graphs with nontrivial first homology).
For a detailed discussion of the map $\wb$, see also Section \ref{sec.fti}.

We should point out that for $M=S^3$ (i.e., for \ihs s)
one can construct sufficiently many invariants of \ihs s to show that
$\wb$ is an isomorphism, over $\BQ$, see \cite{LMO}. The same is true
for finite type (i.e., Vassiliev) invariants of links in $S^3$,
over $\BQ$, see \cite{Ko2}.
However, it is at present unknown whether the map \eqref{eq.graded}
is one-to-one (and thus, an isomorphism), over $\BQ$, for all 3-manifolds.

We now discuss a well-known isomorphism \cite{Ih,O2,Dr,HM}, over
$\BQ$,
for a torsion-free abelian group $A$:
\begin{equation}
\Psi_n: \A^t_n(A) \cong_{\BQ} \LLc_{n+1}(A),
\end{equation}
which will help us relate the Johnson homomorphism to the tree-level
part of \fti s of 3-manifolds.
This map is defined as follows: Fix an oriented uni-trivalent tree
$T$
of degree $n$ (thus with $n+2$ legs, i.e., univalent vertices) and let
$c: \text{Leg}(T)\to A$ be a coloring of its legs. Given a leg $l$ of
$T$, $(T,l)$ is a rooted colored tree to which we can associate an
element
$(T,l)$ of  $\LL_{n+1}(A)$. Due to the $IHX$ relation (see Figure
\ref{lie}),
the function
$$
T \to \sum_{l \in \text{Leg}(T)} c(l) \otimes (T,l)
$$
descends to one $\A^t_n(A) \to A \otimes \LL_{n+1}(A)$ so that its
composition
with $A\otimes \LL_{n+1}(A)\to\LL_{n+2}(A)$ vanishes,
thus defining the map $\Psi_n$.
There is a  map $A \otimes \LL_{n+1}(A) \to \A^t_n(A)$
(defined by sending $a \otimes b \in A \otimes \LL_{n+1}(A)$ to the
rooted tree with one root colored by $a$ and $n$ additional legs
colored
by $c$), which shows that $\Psi_n$ is one-to-one; and by counting
ranks
it follows that it is in fact a vector space isomorphism. It is
unknown
to the authors whether $\Psi_n$ is an isomorphism over $\BZ[1/6]$.

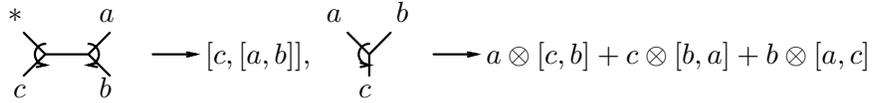
\begin{figure}[htpb]
$$ \printname{lie1}
	\setlength{\unitlength}{0.03\standardunitlength}
	\begin{array}{c}  \hspace{-1.7mm}
         	\raisebox{-8pt}{\begingroup\makeatletter\ifx\SetFigFont\undefined
\def\x#1#2#3#4#5#6#7\relax{\def\x{#1#2#3#4#5#6}}%
\expandafter\x\fmtname xxxxxx\relax \def\y{splain}%
\ifx\x\y   
\gdef\SetFigFont#1#2#3{%
  \ifnum #1<17\tiny\else \ifnum #1<20\small\else
  \ifnum #1<24\normalsize\else \ifnum #1<29\large\else
  \ifnum #1<34\Large\else \ifnum #1<41\LARGE\else
     \huge\fi\fi\fi\fi\fi\fi
  \csname #3\endcsname}%
\else
\gdef\SetFigFont#1#2#3{\begingroup
  \count@#1\relax \ifnum 25<\count@\count@25\fi
  \def\x{\endgroup\@setsize\SetFigFont{#2pt}}%
  \expandafter\x
    \csname \romannumeral\the\count@ pt\expandafter\endcsname
    \csname @\romannumeral\the\count@ pt\endcsname
  \csname #3\endcsname}%
\fi
\fi\endgroup
\begin{picture}(2637,1230)(0,-10)
\thicklines
\path(1245.000,513.000)(1125.000,483.000)(1245.000,453.000)
\put(1125.000,633.000){\arc{300.000}{4.7124}{7.8540}}
\path(405.000,453.000)(525.000,483.000)(405.000,513.000)
\put(525.000,633.000){\arc{300.000}{1.5708}{4.7124}}
\path(225,933)(525,633)(1125,633)(1425,933)
\path(225,333)(525,633)
\path(1125,633)(1425,333)
\path(2025,633)(2625,633)
\path(2505.000,603.000)(2625.000,633.000)(2505.000,663.000)
\put(75,33){\makebox(0,0)[lb]{$c$}}
\put(1275,33){\makebox(0,0)[lb]{$b$}}
\put(0,1083){\makebox(0,0)[lb]{$\ast$}}
\put(1275,1083){\makebox(0,0)[lb]{$a$}}
\end{picture} }
         	\hspace{-1.9mm}
	\end{array}
 [c,[a,b]], \text{ } \printname{lie2}
	\setlength{\unitlength}{0.03\standardunitlength}
	\begin{array}{c}  \hspace{-1.7mm}
         	\raisebox{-8pt}{\begingroup\makeatletter\ifx\SetFigFont\undefined
\def\x#1#2#3#4#5#6#7\relax{\def\x{#1#2#3#4#5#6}}%
\expandafter\x\fmtname xxxxxx\relax \def\y{splain}%
\ifx\x\y   
\gdef\SetFigFont#1#2#3{%
  \ifnum #1<17\tiny\else \ifnum #1<20\small\else
  \ifnum #1<24\normalsize\else \ifnum #1<29\large\else
  \ifnum #1<34\Large\else \ifnum #1<41\LARGE\else
     \huge\fi\fi\fi\fi\fi\fi
  \csname #3\endcsname}%
\else
\gdef\SetFigFont#1#2#3{\begingroup
  \count@#1\relax \ifnum 25<\count@\count@25\fi
  \def\x{\endgroup\@setsize\SetFigFont{#2pt}}%
  \expandafter\x
    \csname \romannumeral\the\count@ pt\expandafter\endcsname
    \csname @\romannumeral\the\count@ pt\endcsname
  \csname #3\endcsname}%
\fi
\fi\endgroup
\begin{picture}(2112,1230)(0,-10)
\thicklines
\path(480.000,453.000)(600.000,483.000)(480.000,513.000)
\put(600.000,633.000){\arc{300.000}{1.5708}{4.7124}}
\path(300,933)(600,633)(900,933)
\path(600,633)(600,333)
\path(1500,633)(2100,633)
\path(1980.000,603.000)(2100.000,633.000)(1980.000,663.000)
\put(0,1083){\makebox(0,0)[lb]{$a$}}
\put(975,1083){\makebox(0,0)[lb]{$b$}}
\put(450,33){\makebox(0,0)[lb]{$c$}}
\end{picture} }
         	\hspace{-1.9mm}
	\end{array}

a \otimes [c,b]+ c\otimes[b,a] +
b\otimes [a,c] $$
\caption{On the left, the map from rooted vertex-oriented trees to
the free
Lie algebra;
on the right the map $\Psi_1$.}\lbl{lie}
\end{figure}

It turns out that a skew-symmetric form $\mathfrak{c}:
A \otimes A \to \BQ$ equips
$\A^t(A)$ with the structure of a graded Lie algebra, by defining the
Lie
bracket
\begin{equation}
\lbl{eq.bracket}
[\Ga,\Ga']^{\mathfrak{c}}= \sum_{a,b}
\mathfrak{c}( a, b )( \Ga_a  \mathrm{ glue }
\Ga'_b ),
\end{equation}
  where the sum is over each leg $a$ of $\Ga$ and $b$ of $\Ga'$ and
$\Ga_a  \mathrm{ glue }
\Ga'_b$ is the graph obtained by gluing the legs $a$ and $b$ of $\Ga$
and $\Ga'$ respectively, with the
understanding that the sum over an empty set is zero.
In other words, $[\Ga,\Ga']^{\mathfrak{c}}$ is the sum of all
{\em contractions} of a leg of $\Ga$ with a leg of $\Ga'$.
This Lie bracket is not new, it has been
observed and used by Morita \cite{Mo1} and Kontsevich \cite{Ko1}
on a close relative of $\A^t(A)$, namely
$\Lie(A)\eqdef A^\ast \otimes
\LL(A)$ (which carries a bracket of degree $-1$).

We now explain the Lie bracket on $\A^t(A)$ from the point of view of
\fti s of 3-manifolds.  Fixing a compact surface $\S_{g,1}$ of genus
$g$ with one boundary component, it follows by definition that
$\M(\sti)$ is generated by homology cylinders over $\S_{g,1}$
and  is a  ring with multiplication
$M_1 \ast M_2$ defined by stacking $M_1$ below $M_2$.
Fix a framed oriented link $\mathfrak b$
in $\sti$ that represents a basis of $H_1(\sti;\BZ)$ and consider
the associated onto map $\wb :\A(\S _{g,1}\times I) \to \GY {}(\sti)$  from
\eqref{eq.graded}, which is expected to be
an isomorphism. Thus, $\A(\sti)$ should be equipped with a  ring
structure. This is the content of the following

\begin{proposition}
\lbl{thm.4}
(i) $\mathfrak b$ induces a homomorphism
$\la \cdot,\cdot \ra^{\mathfrak b}: H_1(\S_{g,1},\BZ)\otimes
H_1(\S_{g,1},\BZ)\to
\BZ$
satisfying\footnote{$ \la \cdot, \cdot \ra^{\mathfrak b} $
will often be denoted by $\la \cdot , \cdot  \ra$ if
$\mathfrak b$ is clear from the context.}
$$ \la a,b \ra^{\mathfrak b} - \la b,a \ra^{\mathfrak b} = a\cdot
b,$$
where $\cdot$ is the natural symplectic form on $H_1(\S_{g,1},\BZ)$.
\newline
(ii) $\A(\sti)$ is a  ring  with $\ast$-multiplication
(depending on $\mathfrak{b}$) defined as follows:
for $\Ga,\Ga' \in \A(\sti)$,
$$
\Ga \ast \Ga'=\sum_{l=0}^\infty \la \Ga,  \Ga' \ra_l ,$$
where
$$\la \Ga, \Ga' \ra_l=
(-1)^l \sum_{a,b}
\prod_{i=1}^l \la a_i, b_i \ra^{\mathfrak b}( \Ga_a  \mathrm{ glue }
\Ga'_b )
$$
is the sum over all ordered subsets  $a=(a_1, \dots
, a_l )$ and
$b=(b_1, \dots , b_l )$ of the set of legs of $\Ga$ and $\Ga'$
respectively,
$\Ga_a \mathrm{ glue } \Ga'_b$ is the graph obtained by gluing the
$a_i$-leg of $\Ga$ to the $b_i$-leg of $\Ga'$, for every $i$,
with the understanding that a sum over the empty set is zero
(thus the multiplication $\ast$ is a finite sum). \newline
(iii) $\A^c(\S _{g,1}\times I)$ is a Lie subring of $\A(\sti)$ with
bracket
defined by
$$
[\Ga,\Ga'] = \Ga \ast \Ga' - \Ga' \ast \Ga=
  \sum_{l=1}^\infty  \la \Ga ,\Ga'\ra_l - \la \Ga', \Ga \ra_l.
$$
(iv) Over $\BQ$, there
is an algebra isomorphism $\mathsf{U}(\A^c(\sti))\cong_{\BQ} \A(\sti)$,
where $\mathsf{U}$ is the universal enveloping algebra functor.
\end{proposition}

\begin{remark}
The leading term $\la \cdot , \cdot \ra_1$ of the $\ast$-multiplication
that involves contracting a single leg is independent of $\mathfrak{b}$
(see also part (iii) of Theorem \ref{thm.5}), whereas the subleading terms
$\la \cdot , \cdot \ra_l$ for $l \geq 2$ depend on $\mathfrak{b}$.
This is a common phenomenon in mathematical physics, analogous to the fact
that differential operators such as the Laplacian or the Dirac depend
on a Riemannian metric, but have symbols independent of it.
\end{remark}


\begin{remark}
It is interesting to compare
$\la a,b \ra^{\mathfrak b}$ to the Seifert matrix of a knot.
Both notions depend
on ``linking numbers'' of ``stacked'' curves, i.e.  curves
pushed in a positive or negative direction and the relation of
$\la \cdot , \cdot \ra^{\mathfrak b}$ and the Seifert matrix to the symplectic
structure on $H_1 (\S )$ is the same. The noncommutativity of
the stacking is reflected by the fact that the form $\la \cdot , \cdot
\ra^{\mathfrak b}$ is not symmetric. For a related appearance of this
noncommutativity, see also \cite{AMR1,AMR2}. Over $\BQ$, the Lie algebra
$\A^c(\sti)$
is closely related to a Lie algebra of chord diagrams on $\S$
considered by
Andersen-Mattes-Reshetikhin in relation to deformation
quantization, loc. cit. We will postpone an explanation of this
relation to a subsequent publication.
\end{remark}

The  ring structure on $\A(\sti)$ would be of little interest were
it not
compatible with the one of $\GY {}(\sti)$ and $\A^t(\sti)$;
this is the content of the following

\begin{theorem}
\lbl{thm.5}
(i) The map $\wb : \A(\sti)\to \GY {}(\sti)$
preserves the  ring structure. \newline
(ii) $\A^l(\sti)$ is a Lie ideal of $\A^c(\sti)$. \newline
(iii) The Lie bracket
of the quotient
$\A^c(\S\times I)/\A^l(\sti)\cong \A^t(\sti)$
is equal to $(-1)^{\mathrm{deg}-1}$ times the Lie bracket of Equation
\eqref{eq.bracket}
using the symplectic form
on $H_1(\sti)\cong H_1(\S_{g,1})$. In particular, it is independent
of the basis $\mathfrak b$.
\end{theorem}

From now on, we will work over $\BQ$.
We now show that the Johnson homomorphism
$\tau: \G\T_{g,1} \to \A^t(\sti)$, or rather its signed version
$\otau\eqdef (-1)^{\mathrm{deg}-1}\tau$,
  can be recovered from the Lie algebra structure on $\A^c(\sti)$,
where $\T_{g,1}(n) \subseteq \Ga_{g,1}[n]$ is the subgroup of the Torelli
group generated by $n$-fold commutators, and $\G_n\T_{g,1}$ denotes the
quotient $\T_{g,1}(n)/\T_{g,1}(n+1)$.
Recall the map $\T_{g,1} \to \M(\sti)$ defined by
changing the parametrization of the top part of homology cylinders,
its linear extension $(I \T_{g,1})^n \to \FY {n}(\sti)$, where $I\T_{g,1}$
is the augmentation ideal of the group ring $\BQ \Ga_{g,1}$, and the induced
algebra map $\G\T_{g,1} \to \GY {}(\sti)$.
The theorem below explains the statement that the Johnson
homomorphism is contained in the
tree-level part of a theory of invariants in $\sti$.

\begin{theorem}
\lbl{thm.6} Given a surface $\S_{g,1}$ as above of genus at least $6$,
there exists a map $\Phi: \G\T_{g,1}
\to \A^c(\sti)$ and commutative diagrams of graded Lie
algebras:
$$
\divide\dgARROWLENGTH by2
\begin{diagram}
\node{\G\T_{g,1}}\arrow{se,r}{\otau}\arrow{s,l}{\Phi} \\
\node{\A^c(\sti)}\arrow{e}\node{\A^t(\sti)}
\end{diagram}
\text{ and }
\begin{diagram}
\node{\G\T_{g,1}}\arrow{se}\arrow{s,l}{\Phi} \\
\node{\A^c(\sti)}\arrow{e,t}{\wb}\node{\GY {}(\sti)}
\end{diagram}
$$
where the left horizontal map is the natural projection on the tree
part. In other words, for $\phi \in \T_{g,1}(n)$ we have
$$ \Phi_n(\phi)=
(-1)^{n-1}\tau_n(\phi) + \mathrm{ loops } $$
in $ \A^c(\sti)$.
\end{theorem}

\begin{remark}
\lbl{rem.closed}
For a closed surface $\S_g$ of genus $g$, there is an identical version
of Proposition \ref{thm.4} and Theorem \ref{thm.5} above. As for
Theorem \ref{thm.6}, given
  a closed surface $\S_g$ of genus at least $6$, there exists
commutative diagrams
$$
\divide\dgARROWLENGTH by2
\begin{diagram}
\node{\G\T_{g}}\arrow{se,r}{\otau}\arrow{s,l}{\Phi} \\
\node{\A^c(\stig)/\Ga_{\omega}}\arrow{e}\node{\A^t(\stig)/\Ga_{\omega}}
\end{diagram}
\text{ and }
\begin{diagram}
\node{\G\T_{g}}\arrow{se}\arrow{s,l}{\Phi} \\
\node{\A^c(\stig)/\Ga_{\omega}}\arrow{e,t}{\wb}\node{\GY {}(\sti)}
\end{diagram}
$$
where $\Ga_{\omega}$ is the ideal of $\A(\stig)$ which is generated
by all elements of the form
$\sum_{i} \Ga_{x_i,y_i,a}$ where $\{x_i,y_i\}$ is a standard symplectic
basis of $H_1(\S_g)$, $a \in H_1(\S_g)$ and $\Ga_{b,c,d}$ denote the
degree $1$ graph
$$\underset{d}{\sideset{^b}{^c}\sfY}$$
with counterclockwise orientation.
\end{remark}

It is natural to ask for a statement of the above theorem
involving general (closed) 3-manifolds $M$.
How does one construct elements in $\FY {}(M)$? Given an imbedding
$\iota:\S\hookrightarrow M$ of a closed surface $\S$ and
  $\phi \in \T(n)$, it was shown in \cite{GGP} that $M-M_\phi \in
\FY {n}(M)$.
The following theorem relates the Johnson homomorphism and the
tree-level part
of \fti s of 3-manifolds.

\begin{theorem}
\lbl{thm.7}
With the above assumptions,  we can find $c^{\mathfrak b}_{\phi,n} \in
\A^l_{n}(M)$
so that we have in $\GY n (M)$:
$$   \wb_n(\Psi_n\i \iota_\ast \otau_n(\phi) + c^{\mathfrak b}_{\phi,n}) =
M-M_{\phi}.$$
\end{theorem}

The above theorem should be compared with \cite[Theorem 6.1]{HM},
where
they show that if a string link $L$ has vanishing $\mu$-invariants of
length less than $n$, then the degree $n$ tree-level part of the
Kontsevich
integral of $L$ is given by the degree $n$ $\mu$-invariants of
$L$. Note that these $\mu$-invariants are Massey products on the
closed
3-manifold obtained by $0$-surgery along the closure of the string
link.

\subsection{Plan of the proof}
\lbl{sub.plan}
The paper consists of two, largely independent sections; the reader
could easily skip one of them without any loss of understanding of the
results of the other. Two notions that jointly appear in Sections
\ref{sec.fti} and
\ref{sec.review} are the Johnson homomorphism and the notion of
homology cylinders.

In Section \ref{sec.fti}, we use combinatorial techniques
that are usually grouped under the name of  \fti s (of knotted objects
such as braids, links, string links or 3-manifolds) or graph
cohomology. A key aspect is the introduction of a Lie
algebra $\A(\sti)$ of graphs and its relation to the Johnson
homomorphism,
via Theorems \ref{thm.4}, \ref{thm.5}, \ref{thm.6} and \ref{thm.7}.

In Section
\ref{sec.review}, we use standard techniques from algebraic and
geometric
topology to prove Theorems \ref{thm.1}  concerning
Massey products in general spaces and closed  $3$-manifolds,in particular,
and Theorem \ref{thm.3} which relates the Johnson
homomorphism to Massey products. In addition, we use standard surgery
techniques adapted to homology cylinders to prove the two realization
Theorems
\ref{th.alf} and \ref{thm.2}.

Finally, in Section \ref{sec.que} we pose a set of questions that naturally
arise in our present study.

\section{Finite type invariants of 3-manifolds}
\lbl{sec.fti}

This section concentrates
on the proof of Theorems \ref{thm.4}, \ref{thm.5}, \ref{thm.6} and
\ref{thm.7}.
The techniques that we use are a combination of geometric and
combinatorial arguments.

\begin{proof}(of Proposition \ref{thm.4})
We only explain the first part. Statements (ii) and (iii) are
obvious and (iv)
follows by a theorem of Milnor-Moore \cite{MM} regarding the structure
of cocommutative graded connected Hopf algebras.

First we arrange that the components of $\eb$ project to immersions in
$\S$ with transverse self-intersections and so that the framing has its
first componenet vector field pointing in the $I$ direction. We call such a
link {\em generic }. Given a two-component
sublink $\{b_1, b_2 \}$ of the framed oriented  link
$\mathfrak b$ in $ \sti$, let
$\{ p(b_1), p(b_2) \}$ denote its projection on $\S\times 0$. Then
$p(b_1)$ and $p(b_2)$ intersect transversely at double points.
Define $\la b_1, b_2 \ra^{\mathfrak b}$ to be the sum with signs over
all
points in $p(b_1) \cap p(b_2)$
that $p(b_1)$ overcrosses $p(b_2)$, according to the convention
$$
\printname{crossing}
	\setlength{\unitlength}{0.03\standardunitlength}
	\begin{array}{c}  \hspace{-1.7mm}
         	\raisebox{-8pt}{
\begingroup\makeatletter\ifx\SetFigFont\undefined
\def\x#1#2#3#4#5#6#7\relax{\def\x{#1#2#3#4#5#6}}%
\expandafter\x\fmtname xxxxxx\relax \def\y{splain}%
\ifx\x\y   
\gdef\SetFigFont#1#2#3{%
  \ifnum #1<17\tiny\else \ifnum #1<20\small\else
  \ifnum #1<24\normalsize\else \ifnum #1<29\large\else
  \ifnum #1<34\Large\else \ifnum #1<41\LARGE\else
     \huge\fi\fi\fi\fi\fi\fi
  \csname #3\endcsname}%
\else
\gdef\SetFigFont#1#2#3{\begingroup
  \count@#1\relax \ifnum 25<\count@\count@25\fi
  \def\x{\endgroup\@setsize\SetFigFont{#2pt}}%
  \expandafter\x
    \csname \romannumeral\the\count@ pt\expandafter\endcsname
    \csname @\romannumeral\the\count@ pt\endcsname
  \csname #3\endcsname}%
\fi
\fi\endgroup
\begin{picture}(3025,1410)(0,-10)
\thicklines
\path(912,483)(537,858)
\path(2488,1007)(2113,1382)
\path(387,1008)(12,1383)
\path(118.066,1319.360)(12.000,1383.000)(75.640,1276.934)
\path(12,483)(912,1383)
\path(848.360,1276.934)(912.000,1383.000)(805.934,1319.360)
\path(2906.934,545.640)(3013.000,482.000)(2949.360,588.066)
\path(3013,482)(2638,857)
\path(2113,482)(3013,1382)
\path(2949.360,1275.934)(3013.000,1382.000)(2906.934,1318.360)
\put(237,33){\makebox(0,0)[lb]{$+1$}}
\put(2337,33){\makebox(0,0)[lb]{$-1$}}
\end{picture}
 }
         	\hspace{-1.9mm}
	\end{array}
.
$$
If $b_1 =b_2$, then we define $\la b_1 ,b_1\ra^{\eb}$ by counting
the self-intersections of $b_1$ in the above manner.
Since $\mathfrak b$ is a basis of $H_1(\sti,\BZ) \cong H_1(\S,\BZ)$,
this defines, by linearity, the desired map
$\la \cdot, \cdot \ra^{\mathfrak b}$.

Since $p(b_1) \cdot p(b_2)$ is the sum with signs over all points
$p(b_1) \cap p(b_2)$, it follows that $\la a,b \ra - \la b,a \ra= a
\cdot b$
for $a,b \in H_1(\S,\BZ)$.
\end{proof}

Before we proceed with the proof of Theorems \ref{thm.5} and
\ref{thm.7},
we need to recall the definition of the map $\wb: \A(M)\to\GY {}(M)$:
Given a colored graph $\Ga$, we will construct an imbedding
of it in $M$ in two steps.

First, we imbed the leaves, as follows. Given the decoration
$x \in H_1(M,\BZ)$, consider its projection $x^{tf} \in
H_1(M,\BZ)/\torsion$
and write $x^{tf}=\sum_{b \in \mathfrak b} n_b [b]$, for integers
$n_b$.
Consider the oriented link $L_x$
obtained by the union (over $b$) of $n_b$ parallel
copies of $b$, where parallel copies of a component of $\mathfrak{b}$
are obtained by pushing off using the framing of $\mathfrak{b}$.
Choose a basing for
$\mathfrak b$, i.e., a set or meridians on each component of
$\mathfrak b$
together with a path to a base point. Join the components of $L_x$
using this basing to construct a knot $K_x$ in $M$.
Apply this construction to every leaf of $\Ga$. Of course, the resulting
link depends on the above choices of basing and joining.

Second, imbed the edges of $\Ga$ arbitrarily in $M$.

This defines an imbedding of $\Ga$ in $M$ (which we denote by the
same name) that also depends on the above choices; however
the associated element $[M,\Ga]$  in $ \GY {}(M)$ ,  where
$$[M,\Ga]=\sum_{\Ga' \subseteq \Ga} (-1)^{|\Ga'|} M_{\Ga'}
$$
is the alternating sum over all $Y$-sublinks $\Ga'$ of $\Ga$, is
well-defined, depending only on the framed link $\eb$.
  This follows
from the following equalities in $\GY {}(M)$ (for detailed proofs
see \cite{Gu2,Ha} and also \cite{GGP}):
\begin{equation}
\lbl{eq.r1}
[ \psdraw{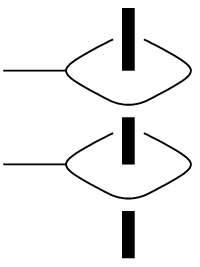}{0.4in}]=[\psdraw{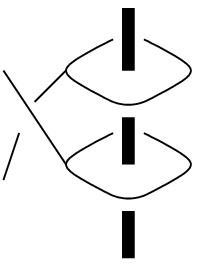}{0.4in}]
\end{equation}
\begin{equation}
\lbl{eq.r2}
  [ \psdraw{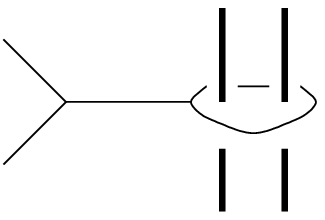}{0.8in}]=[\psdraw{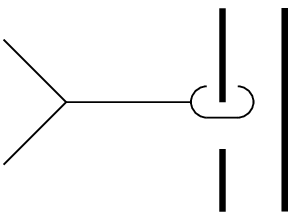}{0.8in}]+[\psdraw{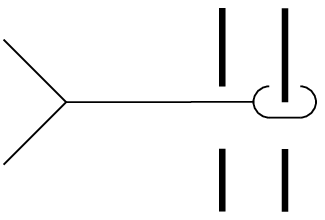}{0.8in}]
\end{equation}
\begin{equation}
\lbl{eq.r3}
  [ \psdraw{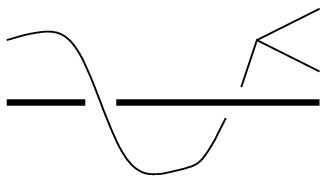}{0.8in}]=[\psdraw{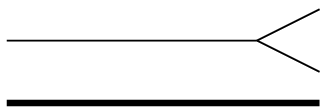}{0.8in}]
\end{equation}
where in the above equalities $[M,\Ga]$ is abbreviated by $[\Ga]$.
Using further identities
in $\GY {}(M)$, one can show loc.cit.
that the map $\Ga \to [M,\Ga]$ factors through
further relations to define a map $\wb: \A(M)\to\GY {}(M)$.

\begin{proof}(of Theorem \ref{thm.5})
For the first part of the theorem, we  begin by
choosing
$\mathfrak b$ in $\sti$ to be a generic link. Note that $\mathfrak b$ can be
recovered
from its projection $p(\mathfrak b)$ together with a knowledge of the signs
at each overcrossing.

Given $\Ga \in \A(\sti)$, let $L^{\mathfrak b}(\Ga)$ be an associated
$Y$-link in $\sti$  such that $W^{\mathfrak b}(\Ga)=
[\sti, L^{\mathfrak b}(\Ga)]$.
Without loss of generality, we can assume that the leaves $l^{\mathfrak
b}(\Ga)$
of $L^{\mathfrak b}(\Ga)$ form a generic link ,
and by abuse of notation, we can write that
$W^{\mathfrak b}(\Ga)=
[\sti, p(L^{\mathfrak b}(\Ga))]$,
with the understanding that we have fixed the signs on the
overcrossings of
$p(L^{\mathfrak b}(\Ga))$.

Now,
given $\Ga,\Ga' \in \A(\sti)$,
let $L^{\mathfrak b}(\Ga)$ and
$L^{\mathfrak b}(\Ga')$ be the associated $Y$-links in $\S\times [0,
1/2]$
and $\S\times [1/2,1]$ respectively, and let $p_i: \sti \to \S\times
\{i\}$
denote the canonical projection.
Then we have that
\begin{eqnarray*}
\wb(\Ga) \cdot \wb(\Ga') & = & [\sti, p_0(l^{\mathfrak b}(\Ga))
\cup p_{1/2}(l^{\mathfrak b}(\Ga'))] \\
& = &  [(\sti )_C, p_0(l^{\mathfrak b}(\Ga))
\cup p_{0}(l^{\mathfrak b}(\Ga'))] \\
& = &  [(\sti )_C, L^{\mathfrak b}(\Ga)
\cup L^{\mathfrak b}(\Ga')]
\end{eqnarray*}
where $C$ is a unit-framed trivial link in $\sti$,   whose
components encircle some  crossings of  $p_0(l^{\mathfrak b}(\Ga))$
and $p_{0}(l^{\mathfrak b}(\Ga'))$, so that surgery on $C$ brings
$p_0(l^{\mathfrak b}(\Ga))$ below $p_0(l^{\mathfrak b}(\Ga'))$.
The following lemma implies that the result of changing an
overcrossing
of $L^{\mathfrak b}(\Ga)$ over $L^{\mathfrak b}(\Ga')$ to an
undecrossing
can be achieved as a difference of the disjoint union of two graphs
minus the disjoint union of two graphs with a leg glued. Together
with Equation \eqref{eq.r2} and the definition of the multiplication
$\Ga \ast \Ga'$, it implies that
\begin{eqnarray*}
[(\sti )_C, L^{\mathfrak b}(\Ga)
\cup L^{\mathfrak b}(\Ga')] & = & \wb(\Ga \ast \Ga')
\end{eqnarray*}
which finishes the proof of the first part of Theorem \ref{thm.5}.

\begin{lemma}
\lbl{lem.2things}
The following identity holds in $\GY {}(\sti)$:
$$ [ \psdraw{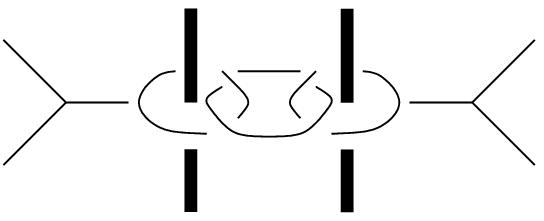}{0.8in}]=[\psdraw{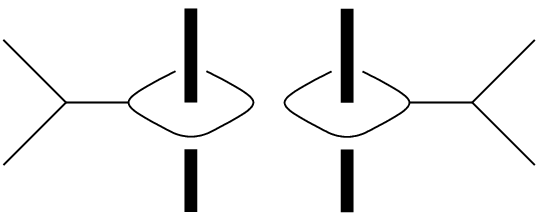}{0.8in}]-[\psdraw{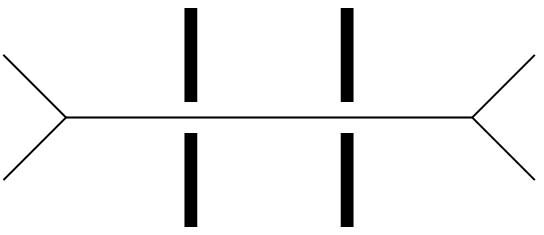}{0.8in}] $$
where the framing of the unknot on the left hand side of the equation
is $+1$ and where we alternate with
respect to the $Y$-links of the figure. The vertical
arcs are arbitrary tubes.
\end{lemma}
\begin{proof}
This follows from the second equality above in $\GY {}(\sti)$
and from
$$ [ \psdraw{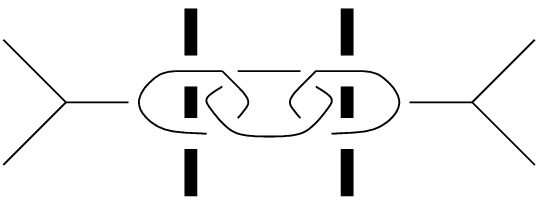}{0.8in}]=-[\psdraw{Y3}{0.8in}] $$
$$ 2 [ \psdraw{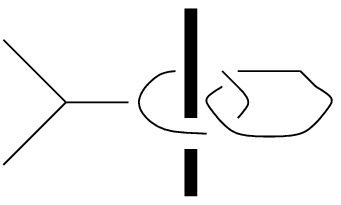}{0.6in}]=0, $$
see \cite{GGP}.
\end{proof}

The second part of Theorem \ref{thm.5} is obvious from the definition
of the Lie bracket.

For the third part, notice that the Lie algebra structure on
the quotient $\A^t(\sti)\cong \A^c(\sti)/\A^l(\sti)
$
is given by
$$
[\Ga,\Ga']=
\la \Ga , \Ga'\ra_1 -\la \Ga , \Ga' \ra_1
=-[\Ga , \Ga']^{\omega},
$$
where the last equality follows from the first part of
Proposition  \ref{thm.4} and where $\omega$ is the symplectic form.
\end{proof}

\begin{proof}(of Theorem \ref{thm.6} and Remark \ref{rem.closed})
For a surface $\S_{g,1}$ of genus $g \geq 3$ with one boundary component,
Johnson \cite{Jo1} introduced a homomorphism $\tau_1: \G_1\T_{g,1} \to
\Lambda^3(H)$, where, following the conventions of
\cite[Chapter 4]{Jo1}, $H=H_1(\S_{g,1},\BZ)$ and
$\La^m(H)$ is identified
with a submodule of the $m$-th tensor power $\mathsf{T}^k(H)$
by defining
$$a_1 \we \cdots \we a_m =\sum_{\pi \in \text{Sym}_m} a_{\pi(1)} \otimes
\cdots\otimes a_{\pi(m)}.$$
In subsequent work, Johnson showed that modulo $2$-torsion his
homomorphism coincides with the abelianization of $\T_{g,1}$, thus
one gets, over $\BQ$,
an onto map of Lie algebras $\LL(\Lambda^3(H))\to \G\T_{g,1}$.

For the rest of the proof we will work over $\BQ$. In \cite[Section 11]{Hain}
Hain proved that for genus $g \geq 6$,
the above map of Lie algebras has kernel generated by quadratic
relations $R_{g,1}$ (Hain's notation for $\G_n\T_{g,1}$ is
$\mathfrak{t}_{g,1}(n)$). Combining the proof of \cite{Hain}
Proposition 10.3 with Theorem 11.1 and Proposition 11.4], it follows that
the relation set $R_{g,1}$ is the symplectic submodule
of $\LL_2(\Lambda^3(H))=\Lambda^2(\Lambda^3(H))$
generated by
$$[ x_1 \wedge x_2 \wedge y_2, x_3 \wedge x_4 \wedge y_4]=0$$
in terms of a standard symplectic basis $\{x_i,y_i\}$ of $H$.

Using the isomorphism $\Lambda^3(H)\cong \A^t_1(\sti)$
given by mapping $a \wedge b \wedge c \in \Lambda^3(H)$
to the degree $1$ graph $\Ga_{a,b,c}$ as in remark \ref{rem.closed},
we obtain a map of Lie algebras $\LL(\Lambda^3(H))\to \A^c(\sti)$.

Since for every
choice of $a \in \{x_1,x_2,y_2 \}$ and
$b \in \{ x_3,x_4,y_4 \}$ we have $a \cdot b =0$, the
first part of Proposition \ref{thm.4} implies
for every basis $\mathfrak b$ we have
$\la a,b \ra^{\mathfrak b} = \la b,a \ra^{\mathfrak b}$. This implies,
by definition, that $[\Ga_{x_1,x_2,y_2}, \Ga_{x_3,x_4,y_4}]=0 \in
\A^c_2(\sti)$, thus obtaining the desired map $\Phi: \G\T_{g,1}
\cong \LL(\Lambda^3(H))/(R_{g,1}) \to \A^c(\sti)$.

Since $\G\T_{g,1}$ is generated by its elements of degree $1$, the
commutativity of the two diagrams follows by their commutativity in
degree $1$; the later follows by definition
for the first diagram, and by the fact that surgery on a $Y$-link
of degree $1$ with counterclockwise orientation and leaves decorated
by $a,b,c$ is equivalent to cutting, twisting and gluing by an element
of the Torelli group (of a surface of genus $3$ with one boundary component,
imbedded in $\S_{g,1}$)
whose image under the Johnson homomorphism is equal to $a \wedge b \wedge c$,
see \cite{GGP}.
This concludes the proof of Theorem \ref{thm.6}.

We now prove the statements in Remark \ref{rem.closed}.
For a closed surface $\S_g$ of genus $g \geq 3$, 
Johnson \cite{Jo1} gave a version of his homomorphism $\tau_1:
\G\T_{g}\to \Lambda^3_0(H)$ where $\Lambda^3_0(H)$ is defined to be
the cokernel of the homomorphism $H\to\La^3(H)$
that sends $x$ to $\omega \we x$, where
$\omega=\sum_{i} x_i \wedge y_i$ is the symplectic form of
$\S_g$, for a choice of symplectic basis.
Working, from now on, over $\BQ$,
Johnson \cite[Chapter 4]{Jo1} gave an identification of $\La^3_0(H)$
with $\ker(C)$, where $C:\La^3(H)\to H$ is given by
$$
C(x \we y \we z)= 2((x\cdot y)z + (y \cdot z)x + (z \cdot x)y)
$$
and $\cdot$ denotes the symplectic form.
Explicitly, we will think of $\La^3_0(H)$ as the submodule of $\La^3(H)$
which is generated by elements of the form  $2(g-1)a -\omega \we
C(a)$ for $a \in \La^3(H)$.
In \cite[Theorems 1.1 and 10.1]{Hain} Hain proved that for genus $g \geq 6$,
there is an isomorphism of graded Lie algebras $\LL(\Lambda^3_0 (H))/(R_g)
\to \G\T_g$ which, in degree $1$, is the inverse of the Johnson
homomorphism, where
  $R_g$ is the symplectic submodule
of $\LL_2(\Lambda^3_0(H))=\Lambda^2(\Lambda^3_0(H))$
generated by the relations
$$[ 2(g-1)x_1 \wedge x_2 \wedge y_2 -x_1 \wedge \omega,
2(g-1)x_3 \wedge x_4 \wedge y_4 -x_3 \wedge \omega ]=0.$$
The slight difference of $2$ in the relations that Hain gave and the ones
mentioned above are due to the difference in the normalization of the
$\we$-product between Hain and Johnson. The restriction of the map
$\LL(\Lambda^3(H))\to \A^c(\sti)\cong \A^c(\stig) \to
\A^c(\stig)/(\Gamma_\omega)$ to $\LL(\Lambda^3_0(H))$ gives a map
$\LL(\Lambda^3_0(H)) \to \A^c(\stig)/(\Gamma_\omega)$ which sends
the relations $R_{g}$ to zero (this really follows from the calculation
of the surface with one boundary component together with the fact that
$a \we \omega \in \Lambda^3(H)$ is sent into the ideal $\Gamma_\omega$
of $A^t_1(\stig)$), thus inducing the desired map $\Phi:
\G\T_{g}\to \A^c(\stig)/(\Gamma_\omega)$.
We claim that $W^{\mathfrak b}(M): \A(\stig) \to \GY {}(\stig)$
maps $\Gamma_\omega$ to zero. This follows from the identity
$$ (\tau_\partial)^{2g-2}=\prod_{i} [ \tau_{x_i}, \tau_{y_i}]$$
of Dehn twists on the mapping class group of $\S_{g,1}$
\cite[Theorem 5.3]{Mo},  where $x_i ,y_i$ refer to the standard
meridian, longitude pairs associated with a symplectic basis of $H_1
(\S_{g,1})$ and $\partial$ is the boundary curve of $\S_{g,1}$. Thus we
have the relation
$$ 1=\prod_i \tau_a [ \tau_{x_i}, \tau_{y_i}] \tau_{a'}^{-1}$$
on the mapping class group of $\S_g$ (where $a,a'$ are simple closed curves
in $\S_{g,1}$ with  isotopic images in $\S_g$), together with the fact
surgery along the $Y$-link $\Ga_{x_i,y_i,a}$ corresponds to the Dehn
twist $\tau_a  [ \tau_{x_i}, \tau_{y_i}] \tau_{a'}^{-1}$ in $\T_g$,
\cite{GGP}.

Since $\G\T_{g}$ is generated by its elements of degree $1$, the
commutativity of the two diagrams follows by their commutativity in
degree $1$; this is shown in the same way as for a surface with one
boundary component.
This concludes the proof of Remark \ref{rem.closed}.
\end{proof}

\begin{proof}(of Theorem \ref{thm.7})
For a closed surface $\S$ of genus at least $6$, Theorem \ref{thm.7} follows
from Remark \ref{rem.closed}
and the following Lemma \ref{lem.nm}, perhaps of independent interest.
For a closed surface $\S$ of genus less than $6$, fix a disk and consider an
imbedding of its complement to a surface $\S'$ in $M$ of genus at
least $6$.
Choose a lifting of $\phi$ to a diffeomorphism of the  punctured
surface that preserves the boundary and extend it trivially to a
diffeomorphism
of $\phi'$ of $\S'$. Since the Johnson homomorphism is stable with
respect
to increase in genus, and since $M_{\phi}=M_{\phi'}$, the result
follows from the previous case.
\end{proof}

\begin{lemma}
\lbl{lem.surface}
If $\Ga \subseteq M$ is an imbedded graph in $M$ with a
distinguished leaf that bounds a surface disjoint from the other leaves
of $\Ga$, then $[M,\Ga]=0 \in \GY {}(M)$.
\end{lemma}

\begin{proof}
First of all, recall that $[M,\Gamma ]=0$ if any leaf bounds a disk
disjoint from the other leaves of $\Gamma$. As explained in \cite{GGP}, an
alternative way of writing Equation
\eqref{eq.r2} is as follows:
\begin{equation}
\lbl{eq.r4}
[ \psdraw{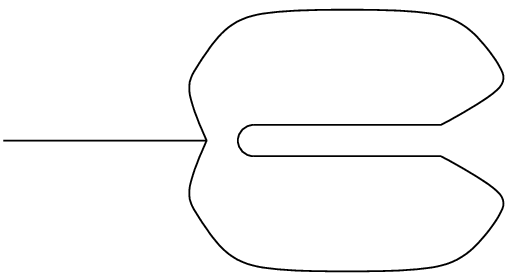}{0.4in}]=[\psdraw{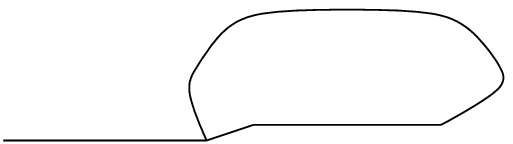}{0.4in}] +
[\psdraw{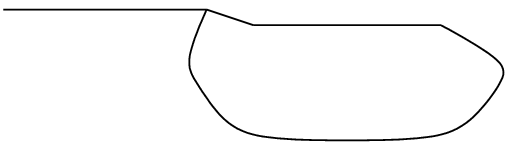}{0.4in}]
\end{equation}
for arbitrary disjoint imbeddings of two based oriented knots in $M$.
Given a based knot $\a$ in $M$, let $\overline{\a}$ denote the based knot
obtained by a push-off of $\a$ in its normal direction (any will do)
followed by reversing the orientation.

The above identity implies that $[M,\Ga_{\overline{\a}}]=-[M,\Ga_{\a}]$
in $\GY {}(M)$, where $\Ga_{\kappa}$ is any imbedded graph in $M$
with a distinguished leaf the based oriented knot $\kappa$ in $M$.

Given $\Ga$ as in the statement of the lemma, it follows that its
distinguished leaf is the connected sum of disjoint based knots of the form
$\a_i \sharp \b_i \sharp \overline{\a_i} \sharp \overline{\b_i}$; thus
it follows from the above discussion that $[M,\Ga]=0$ in $\GY {}(M)$.
\end{proof}

\begin{lemma}
\lbl{lem.nm}
Given an imbedding $\iota: N\to M$ of (not-necessarily closed)
3-manifolds,
and links $\mathfrak{b}(N)$ (resp. $\mathfrak{b}(M$)) in $N$ (resp. $M$)
representing a basis of $H_1(N)$ (resp. $H_1(M)$),
there is an induced map $ \iota_\ast :\A(N) \to \A(M)$
induced by $\iota_*: H_1(N)\to H_1(M)$ on the colorings of
the legs of the graphs and a diagram
$$\begin{diagram}
\node{\A(N)}\arrow{s,l}{W^{\mathfrak{b}(N)}}\arrow{e,t}{\iota_\ast}\node{\A(M)}
\arrow{s,r}{W^{\mathfrak{b}(M)}}\\
\node{\GY {}(N)}\arrow{e,t}{\iota_\ast}\node{\GY {}(M)}
\end{diagram}
$$
that commutes up to $W^{\mathfrak{b}(M)}(\A^l(M))$.
In particular, for $\Ga \in \A^t(N)$ we have
$$
\iota_\ast W^{\mathfrak{b}(N)}(\Ga)=W^{\mathfrak{b}(M)}(\iota_\ast(\Ga))
+ \mathrm{loops}.$$
\end{lemma}

\begin{proof}
Fix a graph $\Ga \in \A^t(N)$. Let $\Ga_L$ (resp. $\Ga_{L'}$) denote
the $Y$-links in $M$ (with leaves $L$ (resp. $L'$)) such that
$$
\iota_\ast W^{\mathfrak{b}(N)}(\Ga)=[M,\Ga_L] \text{ and }
W^{\mathfrak{b}(M)}(\iota_\ast \Ga)=[M,\Ga_{L'}].
$$
It follows by definition that $L$ and $L'$ are homologous
links in $M$. After choosing a common base point for each pair $(L_i,L'_i)$
of components of $L$ and $L'$, it follows that the connected
sum $L_i \sharp \overline{L'_i}$ is nullhomologous in $M$ and thus bounds
a surface $\S_i$ in $M$. The surface $\S_i$ might intersect the other
components of $L$ or $L'$ at finitely many points; however by deleting disks
around the points of intersection of $\S_i$ with $L \cup L'$, we can find
a nullhomotopic based link $L''_i$ and a surface $\S_i'$ disjoint from
$L \cup L'$ with based boundary
such that $L_i \sharp \overline{L'_i}= L_i'' \sharp \partial \S_i'$.
Equation \eqref{eq.r4} and Lemma \ref{lem.surface} imply that
$[M,\Ga_L]-[M,\Ga_{L'}]$ is a sum of terms over $Y$-links $\Ga_\kappa$
in $M$ which are trees, with at least one
component $\kappa$ being nullhomotopic. By choosing a  sequence of
crossing changes (represented by a unit-framed trivial link
$C(\kappa)$) that trivialize $\kappa$ and using Lemma \ref{lem.2things},
it follows that $[M,\Ga_\kappa]=[M,\Ga_{\text{trivial}}]=0$
modulo terms that involve joining
some legs of $\Ga$ (thus modulo terms that involve graphs with
loops), which concludes the proof.
\end{proof}

\section{Massey products and the Johnson homomorphism}
\lbl{sec.review}

\subsection{Universal Massey products}
\lbl{sub.universal}

In this section, homology will be with integer coefficients, unless
otherwise
stated. A useful tool in the proof of Theorem \ref{thm.1} is a
{\em five-term} exact sequence of Stallings \cite{St}:
given a short exact sequence of groups $1 \to H \to G \to K \to 1$,
there is an associated five-term exact sequence
$$
H_2(G) \to H_2(K) \to H/[G,H] \to H_1(G) \to H_1(K) \to 1.
$$
Applying the five-term sequence to the exact sequence
$1 \to R \to F \to G \to 1$, (where $F$ is a free group) we get Hopf's
theorem \cite{Ho}
$$H_2(G)\cong (R \cap [F,F])/[R,F], \text{ and in particular, }
H_2(F/F_{n}) \cong F_{n}/F_{n+1}.$$

In the rest of this section, we will give a proof of Theorem
\ref{thm.1}
and its corollaries.
We will follow a rather traditional notation involving local
coordinates,
\cite{Ma,FS,Dw}.
Let $F$ denote the free group with basis $( x_1 ,\dots, x_m )$; we
will
denote by the same name the corresponding basis of $H=H_1 (F)$.
Let $(u_1 ,\dots, u_m)$ be the dual basis of $H^1 (F)\iso H^1 (F/F_n
)$. The
graded vector space $\bigoplus_n F_n  /F_{n+1} $ has the structure of
a Lie
algebra induced by commutator, and is naturally identified with the
free
Lie algebra $\LL(H)=\bigoplus_n \LL_n (H)$.

We consider the {\em Magnus} expansion, \cite{MKS}.
Let $\pt$ denote the power series ring in
non-commuting variable $\{t_1 ,\dots, t_m\}$. Define $\delta
:F\to\pt$ to be the
multiplicative map defined by $\delta (x_i )=1+t_i$. This is an
imbedding and
induces imbeddings $\delta_n
:\LL_n (H)=F_n  /F_{n+1} \to\pt_n$, where $\pt_n$ is
the subspace of homogeneous polynomials of degree $n$.

We also recall the isomorphism $H_2 (F/F_n  )\iso F_n  /F_{n+1}
=\LL_n (H)$
from above. We
now describe the Massey product structure on $H^ 1 (F/F_n  )$.
For a sequence $I=(i_1,\dots,i_r)$ of numbers $i_j \in \{1,\dots,m\}$
(of length $|I| \eqdef r$), we will let $u_I$ denote the sequence
$(u_{i_1},\dots,u_{i_r})$  and, if each $u_i\in H^1 (F/F_n )$ we let
$\la u_I \ra$ denote the length $r$
Massey
product; we also let  $t_I$ denote the element $\prod_{j=1}^r t_{i_j}$.

\begin{proposition}\label{th.um}
Any Massey product $\la u_I \ra$ of $F/F_n$ vanishes  if
$|I|<n$. The action of any Massey product  $\la u_I \ra$ on
$H_2 (F/F_n  )\iso F_n  /F_{n+1} =\LL_n (H)\sub\pt$ is
determined by the formula
\begin{equation}
\lbl{eq.mp}
  \la u_I\ra \cdot\
t_J =
\begin{cases} 1 &\text{ if }
I=J \text{ and } |I|=n \\ 0 & \text{ otherwise }
\end{cases}
\end{equation}
In other words, the set $\{ \la u_I\ra |  |I|=n \}$ defines the
basis of
$\pt^{\ast}_n$ dual to $\{ t_I |  |I|=n \}$.
\end{proposition}

See \cite{O2} for a slightly less explicit version of this theorem.

\begin{proof}
This follows easily from \cite{FS}. Suppose $w \in F_n$ is some
$n$-fold commutator. Then consider the one-relator group $G=F/\la w
\ra $
and the
projection $p:F/F_n\to G$. Consider $I$ of length $r \leq n$.
By induction we can assume that
$\la u_I  \ra $ is uniquely defined in $F/F_n$ and by \cite{FS},
it is well-defined in $G$. Moreover, by naturality under $p^{*}$,
they take
the same value on $w$. If $r<n$ this is zero by \cite{FS}. Since this
holds
for all $w$ it follows that $\la u_I \ra =0$. If $r=n$
Equation \eqref{eq.mp} follows directly from the formula in \cite{FS}.
\end{proof}

\begin{corollary}
\lbl{cor.j1}
Let  $p :F\to\pi$ be a $2$-equivalence.
Then $p$ is an \nequiv\ if and only if all Massey
products in $H_1 (\pi )$ of length less than $n$ vanish.
\end{corollary}
See also \cite[Proposition 6.8]{CGO}.

\begin{proof}
The ``if'' part follows directly from the above proposition.
To prove the
``only if'' part we proceed by induction on $n$.
The inductive step presents us with a map $\pi \to \pi/\pi_{n-1}
\cong F/F_{n-1}  $; consider the diagram
$$
\divide\dgARROWLENGTH by2
\begin{diagram}
\node{}\node{F/F_n}\arrow{s}\\
\node{\pi}\arrow{e}\arrow{ne,..}\node{F/F_{n-1}}
\end{diagram}$$
The obstruction to lifting this map is the pullback of the
characteristic class
in $H^2 (F/F_{n-1} ;  F_{n-1} /F_n  )$ of the central extension
$F_{n-1} /F_n  \to F/F_{n-1}  \to F/F_n $. But
Proposition~\ref{th.um}
implies that  $H^2 (F/F_{n-1} )$ is
generated by Massey products of length $n-1$ and so the pullback is
zero if
and only if all Massey products of length $n-1$ vanish in $H^2 (\pi
)$.
Thus, we can inductively lift the map to $\pi \to F/F_n$, thus to a
map
$\pi/\pi_n\to F/F_n$, which is still a $2$-equivalence.

On the other hand, since $p$ is a $2$-equivalence, it induces an
onto map $F/F_n\to \pi/\pi_n$, which is also a $2$-equivalence.
Composing
with the map $\pi/\pi_n\to F/F_n$, we get an endomorphism of $F/F_n$
which is a $2$-equivalence. Stalling's theorem \cite{St}
implies that this endomorphism of $F/F_n$ is an isomorphism which
implies that the map $F/F_n\to \pi/\pi_n$ is one-to-one and thus
an \nequiv .
\end{proof}

This proves the first assertion of Theorem \ref{thm.1}.

\begin{corollary}\lbl{cor.lam}
\label{cor.fi}
Suppose $p :F\to\pi$ is an \nequiv .
Then there is an exact sequence:
$$ H_2 (\pi )\overset{\hat{p}}\to F_n  /F_{n+1}
\overset{p_{\ast}}\to\pi_n
/\pi_{n+1} \to 0, $$
where $\hat{p}$ is defined by the formula
$$\hat{p} (\a )=\sum_{I}(\la u_I \ra
\cdot\a )\ t_I, $$
where $\a \in H_2(\pi)$, the summation is over $I$ of length $n$,
$\cdot : H^\ast \otimes H_\ast \to \BZ$ is the
evaluation map,  and where
  the right hand side is asserted to lie in $\LL_n (H)=F_n /F_{n+1}
\sub\pt$.
\end{corollary}

\begin{proof}
Apply Stallings five-term exact sequence to the short exact
sequence of groups $1\to\pi_n\to\pi\to\pi /\pi_n\to 1$ to obtain
$$ H_2 (\pi )\to H_2 (\pi /\pi_n )\to \pi_n /\pi_{n+1}\to 1.  $$
Combining this with the map $p$ gives the commutative diagram

$$\begin{diagram}
\divide\dgARROWLENGTH by2
\node{H_2 (\pi )}\arrow{e}\node{H_2 (\pi /\pi_n )}\arrow{e}\node{\pi_n
/\pi_{n+1}}\arrow{e}\node{1} \\
\node[2]{H_2 (F/F_n )}\arrow{n,l}{\iso}\arrow{e,t}{\iso}\node{F_n
/F_{n+1}}\arrow{n,r}{p_{*}}
\end{diagram}$$

This diagram yields the exact sequence of the corollary, where
$\hat{p}$ is
defined as the composition
$$H_2 (\pi )\to H_2 (\pi /\pi_n )\cong H_2 (F/F_n )\cong
F_n /F_{n+1} $$
To prove the formula for $\hat{p}$ first note that, for any $\a\in H_2
(F/F_n )\iso F_n /F_{n+1}=\LL_n (H)\sub\pt_n$ we have
$$\a =\sum_{I}(\la u_I \ra  \cdot\a )\ t_I, $$
as follows directly from Proposition  ~\ref{th.um}.
But now the corollary follows
from the definition of $\hat{p}$ and naturality.
\end{proof}

This proves the second assertion of Theorem \ref{thm.1} for $K(\pi,1)$
spaces.

\begin{remark}
\lbl{rem.dep}
Given two choices $p, p'$ of maps as in Corollary \ref{cor.fi},
we get a commutative diagram:
$$
\divide\dgARROWLENGTH by2
\begin{diagram}
\node{H_2(\pi)}\arrow{e}\arrow{s,=}\node{F_n/F_{n+1}}\arrow{e}\arrow{s}
\node{\pi_n/\pi_{n+1}}
\arrow{e}\arrow{s,=}\node{0}\\
\node{H_2(\pi)}\arrow{e}\node{F_n/F_{n+1}}\arrow{e}\node{\pi_n/\pi_{n+1}}
\arrow{e}\node{0}
\end{diagram}
$$
where the middle map  is the automorphism of $F_n /F_{n+1}\iso H_2
(F/F_n )$ defined by the composition of isomorphisms
$$H_2 (F/F_n )\overset{p_*}\to H_2 (\pi /\pi_n )\overset{p'_*}\gets H_2
(F/F_n ). $$  In particular, if $p'\con p\ \mod F_2$  then the two exact
rows are identical.
\end{remark}

The third assertion of Theorem \ref{thm.1} for $K(\pi,1)$ spaces
follows from the following

\begin{proposition}
\label{lem.cyc}
Let $p: F \to \pi$ be an \nequiv\ and
$\a_1 ,\dots ,\a_{n+1}\in H^1 (\pi )$. Then we have:
$$\a_1\sms \la \a_2 ,\dots ,\a_{n+1}\ra =\la \a_1 ,\dots, a_n
\ra \sms\a_{n+1}.$$
\end{proposition}

\noindent
See also \cite{Kr} for a related result.

\begin{proof}
We will use Dwyer's formulation \cite{Dw} of the Massey products.
Choose
cocycles $a_i$ representing $\a_i$, for $1\le i\le n+1$. Since we are
assuming all Massey products of length less than $n$ are defined and
vanish,
we can choose cochains $a_{ij}$, for $1\le i<j\le n+2$, with the
exception
of the three cases
$$i=1,j=n+1\qquad i=2, j=n+2 \qquad i=1, j=n+2 $$
  so that $a_{i,i+1}=a_i$ and
$\d a_{rs}=\sum_{r<i<s}a_{ri}\sms a_{is}$. For two of the three
exceptional cases the cochains
$$b_{1,n+1}=\sum_{1<i<n+1}a_{1i}\sms a_{i,n+1}\quad\text{ and }\quad
b_{2,n+2}=\sum_{2<i<n+2}a_{2i}\sms a_{i,n+2}$$
the $b_{ij}$ are cocycles but not necessarily coboundaries. In fact
they
represent the Massey products $\la \a_1 ,\dots ,\a_n \ra $ and $\la
\a_2,\dots
,\a_{n+1}\ra $ respectively.

Thus $\la \a_1 ,\dots ,\a_n \ra \sms\a_{n+1}$ is represented by the
cocycle $b_{1,n+1}\sms a_{n+1,n+2}$
and $\a_1\sms \la \a_2 ,\dots ,\a_{n+1}\ra $ is represented by the
cocycle $a_{12}\sms b_{2,n+2}$.

Now consider the cochain
$$ c=\sum_{1<i<r<n+2}a_{1i}\sms a_{ir}\sms a_{r,n+2} $$
By grouping the terms in one way we see that
\begin{align*}
c &=a_{12}\sms (\sum_{2<r<n+2}a_{2r}\sms
a_{r,n+2})+\sum_{2<i<r<n+2}(a_{1i}\sms a_{ir}\sms a_{r,n+2})\\
&=a_{12}\sms b_{2,n+2}+\sum_{2<i<r<n+2}(a_{1i}\sms\d a_{i,n+2})
\end{align*}
Grouping the terms in another way gives
\begin{align*}
c&=\sum_{1<i<r<n+1}(a_{1i}\sms a_{ir}\sms
a_{r,n+2})+(\sum_{1<i<n+1}a_{1i}\sms a_{i,n+1})\sms a_{n+1,n+2}\\
&=\sum_{1<r<n+1}(\d a_{1r}\sms a_{r,n+2})+b_{1,n+1}\sms a_{n+1,n+2}
\end{align*}
Now subtracting these two formulae for $c$ gives
\begin{align*}
  a_{12}\sms b_{2,n+2}-b_{1,n+1}\sms a_{n+1,n+2}&=\sum_{1<r<n+1}\d
a_{1r}\sms a_{r,n+2}-\sum_{2<i<n+2}a_{1i}\sms\d a_{i,n+2}\\
&=\d (\sum_{2<i<n+1}a_{1i}\sms a_{i,n+2})
\end{align*}
since $a_{12}$ and $a_{n+1,n+2}$ are cocycles. Since the left side
of this
equation represents
$$\a_1\sms \la \a_2 ,\dots ,\a_{n+1}\ra -\la \a_1 ,\dots ,\a_n
\ra \sms\a_{n+1}$$
  and the right side is a coboundary the proof is complete.
\end{proof}

This concludes the proof of Theorem \ref{thm.1} for $K(\pi,1)$ spaces.
The general case follows from the fact that the canonical map
$X \to K(\pi,1)$ induces an onto map $H_2(X)\to H_2(\pi)$.
\qed

\begin{proof}(of Corollary \ref{cor.1})
The first part is immediate from Theorem \ref{thm.1},
using the Poincar\'e duality isomorphism $H^\ast=
H_1(M,\BQ)=H^2(M,\BQ)$.
In local coordinates, it implies (see Corollary \ref{cor.fi})
that the Massey product
$$ \mu_n(M)
\in H \otimes \LL_n(H) \cong \text{Hom}(H^1(M), \LL_n(M))
$$
is given by
\begin{equation}
\lbl{eq.lam}
\mu_n(M) ( u) = \sum_{I} [M]\smf
(u \sms \la u_I\ra ) \ t_I,
\end{equation}
where $u \in H^1(M)$, the summation is over $I$ of length $n$,
$\sms$ indicates cup product and $[M]\smf$ indicates cap
product with the fundamental homology class of $M$.
Let $[ \cdot ]: H \otimes \LL(H)\to\LL(H)$ be the Lie algebra
bracket,
defined by $[a\otimes b]
=[a,b]$, for $a \in H, b \in \LL(H)$.
If we regard $\LL(H)\sub\pt$ then $[ \cdot ]$ can be
expressed by the formula
$$[x_i\otimes c]=(1+t_i )c-c(1+t_i )=t_i c-ct_i.$$
Equation \eqref{eq.lam} implies that
$$\mu_n(M)=\sum_{i,I}x_i\otimes ([M]\smf
(u_i\sms \la u_I\ra ) )\ t_I $$
and so
$$[\mu_n(M)]=\sum_{i,I}
[M]\smf (u_i\sms \la u_I\ra ) \ (t_i t_I-t_I t_i ). $$
Thus Corollary \ref{cor.1} follows from the third assertion of
Theorem \ref{thm.1} (or its coordinate version, Proposition
\ref{lem.cyc}).
\end{proof}

\subsection{Realization results}
\lbl{sub.real}

\begin{proof}(of Theorem \ref{th.alf})
	Given an element $\phi\in \mathrm{A}_0 (F/F_{n} )$ we construct maps
	$f^{\pm}:\S_{g,1}\to K(F/F_{n},1)$,
	where $f^{+}_{*}:\pi_{1}(\S_{g,1} )\to F/F_{n}$ corresponds to the
	canonical projection $p :F\to F/F_{n}$ under the identification
	of $\pi_{1}(\S_{g,1} )$ with $F$, and $f_{*}^{-}=\phi\circ p$. Since
	$\phi_{*}(\omega_g )=\omega_g $, we have
         $f^{+}|\bd\S_{g,1}\simeq f^{-}|\bd\S_{g,1}$ and so
         we can
	combine the two maps to define a map
         $f:\hat{\S}_{g,1}\to K(F/F_{n},1)$, where
	$\hat{\S}_{g,1}$ is the double of $\S_{g,1}$. We would like to
extend
	$f$ to a map $\Phi:M\to\Kk$, for some compact orientable
	$3$-manifold with $\bd M=\hat{\S}_{g,1}$ the obstruction to the
	existence of $\Phi$ is the element $\theta\in\Omega_{2}(F/F_{n})\iso
	H_{2}(F/F_{n})$ represented by $f$, where $\Omega_{*}(F/F_{n})$
	are
	the oriented bordism groups of $F/F_{n}$. Since
	$H_{2}(F/F_{n})\not= 0$ we must be careful in our choices to
	assure that $\theta =0$. Redo the construction of $f^{+},f^{-}$
	and $f$ but using $K(F/F_{n+1},1)$ instead of $\Kk$ and using
	a lift of $\phi$ to an automorphism $\bar \phi$ of $F/F_{n+1}$
	instead of $\phi$. Our restriction on $\phi$ assures that $\bar
	\phi(\omega_g )=\omega_g$ and so we obtain
         $\bar f:\hat{\S}_{g,1}\to K(F/F_{n+1},1)$ and
	an obstruction element $\bar\theta\in\Omega_{2}(F/F_{n+1})$. Now
	this element may not be zero, but since the projection map
	$H_{2}(F/F_{n+1})\to H_{2}(F/F_{n})$ is zero, and clearly
	$\bar\theta$ maps to $\theta$, 	we conclude that $\theta =0$.
	Thus $f$ extends to the desired $\Phi:M\to\Kk$.

	Let $i^{\pm}:\S_{g,1}\to\bd M$ be the obvious diffeomorphisms onto
the
	domains of $f^{\pm}$. It is clear that if $M$ were a homology
	cylinder  over $\S^{+}$, then $\si_{n}(M)=\phi$. But this is not
	necessarily true and so we will perform surgery on the map $\Phi$,
	adapting the arguments in \cite{KM} to our situation. See also
\cite[Theorem 1]{Tu} for similar surgery arguments which are used to show
that any finite $3$-dimensional Poincar\'e complex is homology equivalent
to a closed $3$-manifold.

\begin{lemma}
Suppose $\a\in\ker \Phi_{*}:H_{1}(M)\to H$. Then there exists
$\bar\a\in\pi_{1}(M)$ such that $\bar\a\in\ker \Phi_{*}:\pi_{1}(M)\to
F/F_{n}$ and $\bar\a$ represents $\a$.
\end{lemma}

	\begin{proof}
	If $\bar\a\in\pi_{1}(M)$ is any representative of $\a$, then
	$\Phi_{*}(\bar\a )\in F_{2}/F_{n}$. Choose an element
	$\b\in\pi_{1}(\S^{+})_{2}$ so that $\Phi_{*}(\b )=f^{+}_{*}(\b
	)=\bar\a$. Then $\bar\a\b^{-1}\in\ker \Phi_{*}$ and $\bar\a\b^{-1}$
	represents $\a$.
	\end{proof}

	Thus for any $\a\in\ker \Phi_{*}:H_{1}(M)\to H$ we can do surgery on
	a curve representing $\a$ and extend $F$ over the trace of the
	surgery.

	The first step in killing $\ker \Phi_{*}$ will be to kill the
	torsion-free part. Note that $H_{1}(M)\iso H_{1}(\S^{+})\oplus\ker
	\Phi_{*}$, since $\Phi_{*}\circ i^{+}_{*}$ is an isomorphism, and
so,
	under the canonical map $H_{1}(M)\to H_{1}(M,\bd M),\ \ker \Phi_{*}$
	maps onto $H_{1}(M,\bd M)$. Choose an element $\a\in\ker \Phi_{*}$
	which maps to a primitive element of $H_{1}(M,\bd M)$. Now surgery
	on a simple closed curve $C$ representing $\a$ will produce a new
	manifold $M'$ so that, if $\b\in H_{1}(M')$ is the element
	represented by the meridian of $C$, then
	\begin{equation}\lbl{eq.surg}
	H_{1}(M)/\la \a \ra \iso H_{1}(M')/\la \b \ra
	\end{equation}
	(see \cite{KM}). Since $\a$ is primitive in $H_{1}(M,\bd M)$,
	there is a $2$-cycle $z$ in $M$ whose intersection number
	with $C$ is $+1$. Thus the intersection of $z$ with $M'$ is
	a $2$-chain whose boundary is $\b$. So, by Equation
	\eqref{eq.surg}, $H_{1}(M')\iso H_{1}(M)/\la \a \ra $.

	A sequence of such surgeries will kill the torsion-free part of
	$H_{1}(M,\bd M)$. But this implies that $\ker \Phi_{*}$ is torsion
	by the following simple homology argument. Consider the exact
	sequence:

	$$0\to H_{2}(M)\to H_{2}(M,\bd M)\to H_{1}(\bd M)\to H_{1}(M)\to
	H_{1}(M,\bd M)\to 0
	$$
	Since $\rk H_{2}(M)=\rk H_{1}(M,\bd M)=0$ and $H_{1}(\S^{+})$
	imbeds into $H_{1}(M), \rk H_{2}(M,\bd M)=\rk H_{1}(M)\ge 2g$.
	But, since $\rk H_{1}(\bd M)=4g$, we conclude that $\rk
	H_{1}(M)=2g$. Therefore
	$$ 2g=\rk H_{1}(M)=\rk H_{1}(\S^{+})+\rk\ker \Phi_{*} $$
	and so $\rk\ker \Phi_{*}=0$.

	We now follow the argument in \cite{KM} to kill the torsion group
$T=\ker \Phi_*$. The linking pairing $l:T\otimes T\to\BQ /\BZ$ is
non-singular
since $T=\tor H_1 (M)$ maps isomorphically to
$\tor H_1 (M,\bd M)=H_1 (M,\bd M)$. According to \cite[Lemma 6.3]{KM}
if,
for $\a\in T$, $l(\a ,\a )\not=0$, then we can choose the normal
framing to
any closed curve $C$ representing $\a$ so that the element $\b\in H_1
(M')$
is of finite order smaller than the order of $\a$. Thus the torsion
subgroup
of $H_1 (M')$ is smaller than $T$. Continuing in this way we reach
the point
where all the self-linking numbers are $0$. According to \cite[Lemma
6.5]{KM}
this implies that $T$ is a direct sum of copies of $\BZ /2$. Now
choose any
non-zero element $\a\in T$. We will show that surgery on $\a$ reduces
the
rank of $H_1 (M,\bd M;\BZ /2)$. Denote by $V$ the trace of the
surgery and $M'$
  the result of the surgery. Then we have a diagram of homology groups
(coefficients in $\BZ /2$) with exact row:
$$
\divide\dgARROWLENGTH by2
\begin{diagram}
\node{H_2 (V,\bd M')}\arrow{e}\node{H_2 (V,M')}\arrow{e}\node{H_1
(M',\bd M')}\arrow{e}\node{H_1 (V,\bd M')}\arrow{e}\node{0}\\
\node{H_2 (M,\bd M)}\arrow{n}\arrow{ne,..}
\end{diagram}
$$
	Now $H_1 (V,\bd M')\iso H_1 (M,\bd M)/\la \a \ra $ and so has rank
one less than $H_1 (M,\bd M)$. Since $H_2 (V,M')$ is generated by the
transverse disk bounded by the meridian curve representing $\b$, the
dotted
arrow can be interpreted as the ($\BZ /2$) intersection number with
$\a$. By
Poincar\' e duality this map is non-zero and so
$H_1 (M',\bd M')\iso H_1 (V,\bd M')$, proving the claim. As in
\cite{KM}
we can assume the normal framing chosen so that $\b$ has order $2$ or
$\infty$. Thus the possibilities for $H_1 (M',\bd M';\BZ)$ are either
$\BZ\oplus (s-2)\BZ /2$ or $\BZ /4\oplus (s-2)\BZ /2$, where
$s=\rk H_1 (M,\bd M)$. We can then do a surgery to kill the $\BZ$
factor,
in the first case, or reduce the order of $H_1 (M',\bd M';\BZ )$,
in the second case. Continuing this way we eventually kill $\ker
\Phi_*$,
producing the desired $(M,\Phi)$.
\end{proof}

\begin{proof}(of Proposition \ref{prop.DA})
Let $\LLc^{\mathrm{a}}_n(H)$ denote the kernel of the natural
projection
$\mathrm{A}_0(F/F_{n+1})\to\mathrm{A}_0(F/F_{n})$.
We first construct a map $D_n:\LLc^{\mathrm{a}}_n(H)\to \LLc_n(H)$ as
follows.
If $h \in \LLc^{\mathrm{a}}_n(H)$ we can write $h(a)=a \psi(a)$,
where $\psi(a) \in
F_n/F_{n+1}\cong \LL_n(H)$.
Then, we define
$D_n(h)([a])=\psi(a)$,
where $[a] \in H$ and $a \in F/F_n$ is a lift of $[a]$.
Using the isomorphism $\text{Hom}(H, \LL_n(H)) \cong H \otimes
\LL_n(H)$
this defines a map (denoted by the same name)
$$\LLc^{\mathrm{a}}_n(H)\to H \otimes \LL_n(H)$$
with corresponding description in local coordinates given by
$$D_n(h)=\sum_{i} x_i \otimes \psi(y_i)- y_i \otimes \psi(x_i) \in
H \otimes \LL_n(H).$$
If $h\in \LLc^{\mathrm{a}}_n (H)$, as above, then $\prod_i [h(x_i
),h(y_i )]\con\prod_i
[x_i ,y_i ]\mod F_{n+2}$ and so
$$\prod_i [x_i\psi (x_i ),y_i\psi (y_i )]\con\prod_i [x_i ,y_i ][\psi
(x_i
),y_i ][x_i ,\psi (y_i )]\mod F_{n+2} $$
Therefore $\prod_i [\psi (x_i ),y_i ][x_i ,\psi (y_i )]\in F_{n+2}$,
which implies that $D_n(h) \in \LLc_n(H)$.

It is clear that $D_n$ is one-to-one.
We now show that it is onto.
Suppose we have an element $\theta =\sum_i (x_i\otimes\a_i
-y_i\otimes\b_i )\in\LLc_n(H) $. Lift $\a_i ,\b_i$ into $F_n$
(denoted by
the same symbols) and define an endomorphism $h$ of $F$ by
$$h(x_i )=x_i\a_i ,\ h(y_i )=y_i\b_i .$$
It follows by Stalling's theorem \cite{St} that
$h$ induces an automorphism of $F/F_{n+1}$ which restricts to the
identity automorphism of $F/F_n$. We note that
$$h(\prod_i [x_i ,y_i ])=\prod_i [x_i\a_i ,y_i\b_i ]\con\prod_i [x_i
,y_i
][x_i ,\a_i ][\b_i ,y_i ]\mod F_{n+2}$$
But $\prod_i [x_i ,\a_i ][\b_i ,y_i ]$ represents
the image of $theta$ under the Lie bracket
$H \otimes \LL_n(H) \to \LL_{n+1}(H)$, which vanishes since
$\theta \in \LLc_n(H)$; thus $\prod_i [x_i ,\a_i ][\b_i ,y_i ] \in
F_{n+2}$.
This shows that $h\in\mathrm{A}_0 (F/F_{n+1})$ and clearly
$D_n (h)=\theta$.

The fact that the projection
$\mathrm{A}_0(F/F_{n+1})\to\mathrm{A}_0(F/F_{n})$ is onto follows
immediately from Theorem \ref{th.alf}. It is not hard, however, to give a
direct argument; we leave this as an exercise for the reader.
Finally, it is clear by the definitions that the diagram in
Proposition \ref{prop.DA} commutes, and that the sequence below it is exact.
\end{proof}

  \begin{remark}
The action of $\mathrm{A}_0(F/F_{n})$ on $\LLc^{\mathrm{a}}_n(H)\iso
\LLc_n(H)$ induced by conjugation by elements of $ \mathrm{A}_0(F/F_{n+1})$
coincides with the natural action of
$\mathrm{A}_0(F/F_{2})\iso\text{Sp}(2g,\BZ )$ on $ \LLc_n(H)$, via the
projection $ \mathrm{A}_0(F/F_{n})\to \mathrm{A}_0(F/F_{2})$.
\end{remark}

\begin{proof}(of Theorem \ref{thm.2})
Let $M\in\H_{g,1}$ and define $S(M)^o=T_+\cup M\cup T_-$, where
$T_{\pm}$
are two copies of the solid handlebody $T$ of genus $g$, which are
attached
to $\bd M$ via the diffeomorphisms $i^{\pm}$ so that, referring to a
basis
$\{ x_i ,y_i\}$ of $F$ corresponding to a symplectic basis of $H$,
the $\{
x_i\}$ are represented by the boundaries of meridian disks in $T$.
Thus
$\pi_1 (T)=F'$, the free group generated by $\{ y_i\}$ (or, more
precisely,
their images in $\pi_1 (T)$). $S(M)^o$ is a 3-manifold with boundary
$S^2$,
which we can fill-in to obtain a closed 3-manifold $S(M)$.
If $M\in\H_{g,1}[n]$, then the inclusion
$T_+\sub S(M)$ induces an isomorphism $p :F'/F'_n\cong \pi_1
(S(M))/\pi_1
(S(M))_n$ and we can consider $\mu_n(S(M),p) \in \LLc_n(H')$,
where $H'=H_1 (F')$.
Suppose that $\si_{n+1}(M)=h\in\mathrm{A}_0 (F/F_{n+1})$.
Set $a_i =\rho (h(x_i ))\in
F'_n$, where $\rho :F\to F'$is the projection defined by $\rho (x_i
)=1$.
Then $ \mu_n(S(M),p)=\sum_i [y_i ]\otimes [a_i ]$, where $[y_i ]
\in H' ,[a_i ] \in \LL_n(H')$
are the  classes represented by $y_i ,a_i$. This assertion is just the
obvious generalization of Corollary \ref{cor.1} and the proof is the
same.

Now let $\sum_i [y_i ]\otimes [\l_i ]$ be an arbitrary element in
$\LLc_n(H')$, where $\l_i\in F'_n$, i.e. $\prod_i
[y_i ,\l_i ]\in F_{n+2}$. We want to construct $(N,p )$ such that
$\mu_n(N,p)=\sum_i [y_i ]\otimes [\l_i ]$. Consider the
endomorphism $h$ of $F$ defined by
\begin{equation}\lbl{eq.O}
\begin{split}
h(x_i )&=x_i\l_i \\
h(y_i )&=\l_i\i y_i\l_i.
\end{split}
\end{equation}
Denote also by $h$ the induced
automorphism of $F/F_{n+1}$.
To see that $h\in\mathrm{A}_0 (F/F_{n+1})$, we compute
$$\prod_i [h(x_i ),h(y_i )]=\prod_i (x_i y_i x_i^{-1}\l_i^{-1}
y_i^{-1}\l_i
)=\prod_i [x_i ,y_i ][y_i ,\l_i^{-1} ]\con\prod_i [x_i ,y_i ]\mod
F_{n+2}.$$
Therefore, by Theorem \ref{th.alf}, $h=\a_n (M)$ for some
$M\in\H_{g,1}$. Since $\l_i\in F'_n$ we have $h\in \LLc^{\mathrm{a}}_n(H)$
and so  $M\in\H_{g,1}[n]$.
By the discussion above, $\mu_n(N,p)=\sum_i [y_i ]\otimes
[\l_i ]$.

\end{proof}

For completeness, we close this section by a sketch of a  more direct
proof of Theorem \ref{thm.2} using the results of \cite{O1,IO}.
Similar arguments can be found in \cite{CGO}.

\begin{lemma}
\lbl{lem.real}
For every $\a \in H_3(F/F_n)$ there is a closed 3-manifold $M$
and an \nequiv\ $p:F\to\pi\eqdef \pi_1(M)$ such that $p_\ast[M]=\a$.
\end{lemma}

\begin{remark}
Since $\Omega_3 (F/F_n )\iso H_3 (F/F_n )$,  it follows
that every element $\a \in \Omega_3 (F/F_n )$ is represented by
some closed 3-manifold $M$ and map $p :\pi\eqdef \pi_1(M)\to F/F_n$
so that
$p_* [M]=\a$. The point is to arrange that
$p_* :H_1 (M)\iso H_1
(F/F_n )$, which would imply that $p$ is an \nequiv .
\end{remark}

\begin{proof}
We apply the constructions and results of \cite{O1}.  Consider the
mapping cone $K_n$ of the natural
map $K(F,1)\to K(F/F_n,1)$ of Eilenberg-MacLane spaces. $K_n$ is
constructed from $K(F/F_n,1)$ by adjoining
$2$-cells $e^2_i$ along the generators $x_i\in F\twoheadrightarrow
F/F_n$.
Then $K_n$ is simply-connected and $H_i (K_n )\iso H_i (F/F_n )$ for
$i\ge
2$. So we have the Hurewicz epimorphism $\rho :\pi_3 (K_n
)\twoheadrightarrow
H_3 (K_n )\iso H_3 (F/F_n )$, where  $\pi_3$
denotes the third homotopy group. Suppose $\rho (\theta )=\a$. Then
the
Pontrjagin-Thom construction gives us a map $f:S^3 \to K_n$
representing
$\theta$, such that, if $x_i\in e^2_i$ is some interior point, then $f\i
(x_i )=L_i$, a zero-framed imbedded circle in $S^3$, see \cite{O1}.
Now let
$M^3$ be the result of framed surgery on $S^3$ along the $\{ L_i\}$.
Then
$f$ induces a map $p :M\to K(F/F_n )$ and it is clear that $p_*
[M]=\a$. Finally we note that $p_* :H_1 (M)\iso H_1 (F/F_n )$.
\end{proof}

From another viewpoint, we have
defined
a homomorphism $\mu_n ' :H_3 (F/F_n )\to \LLc_n(H)$ by the formula
\begin{equation}\lbl{eq.hom}
\mu_n '(\a )=\sum_{i,I} x_i\otimes (\a\smallfrown (u_i\sms
\la u_I \ra ))t_I
\end{equation}
and so, by Corollaries \ref{cor.1} and \ref{cor.lam} $\mu_n (M,p )=\mu_n '
(p_* [M])$. It follows  from Theorem \ref{thm.1} that the kernel of $\mu_n
'$ is precisely
  the image of
$H_3 (F/F_{n+1})\to H_3 (F/F_n )$, inducing, therefore, an injection
$\cok ( H_3 (F/F_{n+1})\to H_3 (F/F_n ))
\rightarrowtail \LLc_n(H)$. But it is shown in
\cite{O1,IO} that both sides are free
finitely generated abelian groups of the same rank and so they are
isomorphic.  In fact K. Orr has pointed out to us that a
straightforward
examination of the spectral sequence of the group extension
$F_{n}/F_{n+1}\to F/F_{n+1}\to F/F_{n}$ shows that the rank of
$\cok ( H_3 (F/F_{n+1})\to H_3 (F/F_n ))$ is at least that of $\LLc_n(H)$.
\qed

\subsection{Massey products and the Johnson homomorphism}

In this section we will give a proof of Theorem \ref{thm.3}.
We first recall the definition of the Johnson homomorphism. Let
$\mathrm{A}_0 (F/F_n)$ be as in Section \ref{sub.homcyl}.
We define the Johnson homomorphism
$$\tau_n :\mathrm{A}_0 (F/F_n)\to
\hom (H,\LL_n (H))$$
where $H=H_1 (F)$,
  as follows. If $\phi\in \mathrm{A}_0 (F/F_n)$,
then $\tau_n (\phi )\cdot [\a ]=\a
\phi(\a\i )\in F_n /F_{n+1}\iso\LL_n (H)$, where $[\a ]\in H$ is the
homology
class of $\a\in F$. Let $K(\phi )\sub F$ be the normal closure of all
elements of
the form $\a\phi(\a\i )$ for $\a\in F$; note that $K(\phi )\sub F_n$ if
$ \phi\in \mathrm{A}_0 (F/F_n)$.
Let $\pi =F/K( \phi )$. Then the projection $F\to\pi$
induces an isomorphism $p_{ \phi}:F/F_n\overset\iso\to \pi /\pi_n$ and
the
homomorphism
$\hat p_{ \phi }$ defined (as $\hat p$) in Corollary~\ref{cor.fi} is given
by the
composition
$$H_2 (\pi )\iso K( \phi ) /[F,K(\phi ) ]\to F_n /F_{n+1}\iso\LL_n (H) $$
The first isomorphism is given by Hopf's theorem. Now consider the
natural
homomorphism $i:H\to H_2 (\pi )$ defined by $i([\a ])=[\a\phi(\a )\i
]$, then
\begin{equation}\lbl{eq.J}
\tau_n (\phi )=\hat p_{\phi}\circ i.
\end{equation}
Now let $\S$ be a compact orientable surface of genus $g$ with one
boundary
component. Then $\pi_1 (\S )\iso F$, the free group on $2g$
generators.
By Nielsen's theorem $\Ga_{g,1}\sub \mathrm{A}(F)$ and so $\Ga_{g,1}
[n]=\Ga_{g,1}\cap \mathrm{A}_0 (F/F_n)$.

For $ \phi\in\Ga_g$ define an associated closed $3$-manifold $T_{\phi}$ as
follows. Let $T'_{\phi}$ be the mapping torus of $ \phi$, i.e. $I\times\S$
with
$1\times\S$ identified with $0\times\S$ by the homeomorphism $(1,x)\to
(0, \phi(x))$ for every $x\in\S$. Since
$ \phi |\bd\S =\text{id}$ there is a canonical identification of $\bd
T'_{\phi}$ with
the torus
$S^1\times S^1$. $T_{\phi}$ is defined by pasting in $D^2\times S^1$. Note
that $\pi_1 (T_{\phi})\iso F/K(\phi )$, where $K(\phi )$ was defined above.

A theorem of the following sort was first suggested by Johnson
\cite{Jo1} and a
proof (which we sketch, for completeness) was given in \cite{Ki}.

\begin{proposition}\lbl{cor.2}
$\tau_n (\phi )=\mu_n (T_{\phi })$, where $H_1 (\S )$ and $H^1 (T_{\phi})$ are
identified by
the string of isomorphisms
$$H^1 (T_{\phi})\overset\iso\to H^1 (\S )\iso H_1 (\S )$$
induced by the inclusion $\S\sub T_{\phi}$ and Poincar\'e duality for
$\S$.
\end{proposition}

\begin{proof}

It follows from the definitions that we need to establish the
commutativity of
the following diagram

$$\begin{diagram}
\node{H_1 (\S )=H}\arrow[2]{e,t}{i}\arrow{s,l,<>}{\iso}\node[2]{H_2
(\pi )}\\
\node{H^1 (\S )}\node{H^1
(T_{\phi})}\arrow{w,t}{j^*}\arrow{e,t,<>}{\iso}\node{H_2
(T_{\phi})}\arrow{n,r}{\theta}\arrow{wnw,t}{j_*}
\end{diagram}$$

where $i$ is defined in equation~\eqref{eq.J}, $j:\S\sub I \times 0 \sub
  I\times\S\to
T_{\phi}$
is the inclusion map inducing $j^*$ and the Gysin homomorphism $j_*$,
and
$\theta$ is the Hopf map $H_2 (X)\to H_2 (\pi_1 (X))$ when $X=T_{\phi}$.

Now suppose $\a\in\pi_1 (\S )$ is represented by a $1$-cycle $z$ in
$\S$. since $\phi(z)$ is homologous to $z$, there exists a $2$-chain $c$
in $\S$ such that $\bd c=\phi(z)-z$. Consider the chain $I\times z$ in
$I\times \S\to T_{\phi}$ with boundary $1\times z-0\times z$. Since
$0\times
z$ is identified with $1\times \phi(z)$ in $T_{\phi}$, the chain $\xi
=I\times
z+1\times c$ is a $2$-cycle in $T_{\phi}$---let $\b$ denote  its homology
class in $H_2 (T_{\phi})$.

It follows from the definition of $\theta$ that $\theta (\b )=i([\a
])$.
On the other hand $j_* (\b )=[\a ]$ since the $2$-cycle $\xi$
intersects
$\e\times\S$ transversely in the $1$-cycle $\e\times z$, if $0<\e <1$.
This establishes the desired commutativity.
\end{proof}

Combining Proposition~\ref{cor.2} and Corollary~\ref{cor.lam} we have

\begin{corollary}\lbl{cor.mas}
If $\{ x_i\}$ is a basis of $H=H_1 (\S )\iso H^1 (T_{\phi})$ and $\{
u_i\}$
is the dual basis of $H_1 (T_{\phi})$, then
$$ \tau_n (h)=\sum_{i,I}x_i\otimes
([T_{\phi}]\smallfrown
(u_i\smallsmile \la u_{I}\ra )) t_{I} $$
using the
inclusion $\LL_n (H)\sub\mathbb Z[\![t_1 ,\dots,t_{2g}]\!]_n$.
\end{corollary}

We now turn to the proof of Theorem \ref{thm.3}. Fix an \nequiv\
$p: F \to \pi_1(M)$, an imbedding
$\iota: \S \to M$ of a (closed) surface in a (closed) 3-manifold $M$,
and
an element $\phi \in \Ga_g[n]$.
The homomorphisms $f: \pi_1(M)\to F/F_{n+1}$
and $f_{\phi}: \pi_1(M_{\phi})\to F/F_{n+1}$ induced by the \nequiv\
$p$
(in the discussion before the statement of the theorem) determine
maps $f: M \to K(F/F_{n+1},1)$ and $f_{\phi}: M_{\phi} \to
K(F/F_{n+1},1)$
in the Eilenberg-MacLane space $K(F/F_n,1)$. We will construct a
cobordism
$F: V^4 \to K(F/F_{n+1})$ between $f$ and $f_{\phi}$.

First we push $\iota(\S)$ in the positive and negative normal
directions in
$M$ to obtain two copies of $\iota(\S)$,  $\S_p$ and $\S_n$,
respectively.
Then we attach $I\times\iota(\S)$ to $M$ along its boundary by
identifying
$0\times\iota(\S)$ to $\S_n$ by the homeomorphism $(0,x)\to x$, for
$x\in\iota(\S)$, and $1\times\iota(\S)$ to $\S_p$ by $(1,x)\to
\phi(x)$.
We can
now thicken up $I\times\iota(\S)$ and $M$ to obtain a manifold $W$,
see
Figure \ref{W}.

\begin{figure}[htpb]
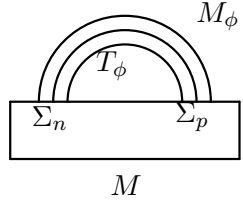

$$ \printname{W}
	\setlength{\unitlength}{0.04\standardunitlength}
	\begin{array}{c}  \hspace{-1.7mm}
         	\raisebox{-8pt}{\input draws/W.tex }
         	\hspace{-1.9mm}
	\end{array}
 $$
\caption{The manifold $W$. The handle shown represents a thickening
of $\iota(\S)$, shown as the core.}\lbl{W}
\end{figure}

Note that the boundary of $W$ consists of three components:
$M$,$M_\phi$ and the mapping torus $T_{\phi}$ of $\phi$ on
$\iota(\S)$. We make one small
modification to obtain $V$. Remove a disk $D$ from $\iota(\S)$ to
obtain
a surface $\S^o$ with boundary and lift $\phi$ to a diffeomorphism
$\phi^o$
of $\S^o$.
Then $S^1\times D$ is naturally
imbedded in $T_\phi$, and so in the boundary of $W$. We attach
$D^2\times D$ to $W$ along $S^1\times D$ to obtain $V$.

The boundary of $V$ is now given by $\partial V=M_\phi
-T_{\phi^o}-M$. It is
not difficult to check that the inclusions $M\to V\gets M_\phi$
are $(n+1)$-equivalences since
  since $\phi\in\Ga_g [n+1]$.
Then $f,f_{\phi}$ extend to a map $F:V\to K(F/F_{n+1} ,1)$.
In addition the inclusions
$\S^o\sub T_{\phi^o}\sub V$ induces isomorphisms
$$F/F_{n+1}\iso\pi_1 (\S^o )/\pi_1 (\S^o)_{n+1} \cong
\pi_1 (T_{\phi^o})/\pi_1 (T_{\phi^o})_{n+1} \cong \pi_1 (V)/\pi_1
(V)_{n+1}$$
  So the bordism $(F,V)$ implies that $(f_{\phi})_* [M_\phi ]=f_* [M]
+(p_{\phi})_*
[T_{\phi^o}] \in \Omega_3(F/F_{n+1}) \cong H_3(F/F_{n+1})$.

Applying the map of Equation \eqref{eq.h3} concludes the proof of
Theorem
\ref{thm.3}.
\qed

\begin{remark}
\lbl{rem.hcc}
Theorem \ref{thm.3} and Proposition \ref{cor.2} can be
generalized easily to the context of homology cylinders, by essentially the
same arguments. For any homology cylinder $N$ there is an obvious notion
of mapping torus $T_N$---then Proposition \ref{cor.2} can be rephrased,
replacing $T_{\phi}$ by $T_N$, $\phi\in\Ga_{g,1}[n]$ by $N\in\H_{g,1}[n]$
and $\tau_n$ by $\varsigma_n$. For Theorem \ref{thm.3} we consider
$N\in\H_{g,1}$ and an imbedding $\S\sub M$, where $\S$ is a closed
orientable surface of genus $g$. We can cut $M$ open along $\S$ and paste
in $\bar N$---where $\bar N$ is obtained from $N$ by filling it in, in the
obvious way, to get a homology cylinder over $\S$---and so obtain a new
manifold $M_N$. Then Theorem \ref{thm.3} can be rephrased, replacing
$M_{\phi}$ by $M_N$.
\end{remark}

\begin{remark}\lbl{rem.fil}
We can consider another filtration of $\H_{g,1}$. Let $\H_{g,1}(n)$ denote
the set of all homology cylinders $n$-equivalent to $\sti$ (see Section
\ref{sub.fti}). Thus if $M\in\H_{g,1}(n)$ then $M-(\sti )\in\F_n^Y (\sti
)$. It is easy to see that $\H_{g,1}(n)$ is a subsemigroup of $\H_{g,1}$.
It is a natural conjecture that the quotient $\H_{g,1}/ \H_{g,1}(n)$ is a
group.

  According to \cite{GGP} we get the same filtration if we ask that $M$ be
obtained from $\sti$ by cutting open along some imbedded closed orientable
surface $\S '\sub\sti$ and reattaching by some element of $\T_n $, the
$n$-th lower central series subgroup of the Torelli group $\T$ of $\S '$.
It is clear that $\H_{g,1}(n)\sub \H_{g,1}[n]$ since the effect of cutting
and reattaching in $M$ by an element of $\T_n$ does not change $\pi_1
(M)/\pi_1 (M)_n$. If $\G\H_{g,1}(*)$ and $\G\H_{g,1}[*]$ denote the
associated graded groups of these filtrations then we have a natural map $
\G\H_{g,1}(*)\to \G\H_{g,1}[*]$. Note also that the natural homomorphism
$\Ga_{g,1}\to\H_{g,1}$ induces maps $(\T_{g,1})_n\to\H_{g,1}(n)$ and so
$\G (\T_{g,1})_*\to\G\H_{g,1}(*)$.

Putting this, and some of the other maps constructed in this paper, all
together, we have a commutative diagram:

$$
\divide\dgARROWLENGTH by2
\begin{diagram}
\divide\dgARROWLENGTH by2
\node[3]{\G (\T_{g,1})_*}\arrow[2]{sw}\arrow[2]{s}\arrow[2]{e}\node[2]{\LLc_*
(H)}\arrow[2]{s,r}{\iso}\\ \\
\node{\A^c (\sti
)}\arrow[2]{ese}\arrow[2]{e,..}\arrow[2]{se}\node[2]{\G\H_{g,1}(*)}
\arrow{s,-}\arrow[2]{e} \node[2]{\G\H_{g,1}[*]}\arrow[2]{s,r}{\iso}\\
\node[3]{}\arrow{s}\\
\node[3]{\G^Y\M (\sti )}\node[2]{\A^t (\sti )}
\end{diagram}
$$
The dotted arrow denotes a conjectured lifting.
\end{remark}

\section{Questions}
\lbl{sec.que}

It is well known that there is a set of moves that generates (string) link
concordance, \cite{Tr}. These moves, together with
the existence of the Kontsevich integral, were the key to the proof
that the tree-level part of the Kontsevich integral of string-links
is given by Milnor's invariants, or equivalently, by Massey products,
see \cite{HM}.

\begin{question}
Is there a set of local moves that generates homology cobordism of homology
cylinders?
\end{question}

\begin{question}
Does Theorem \ref{thm.6} generalize to homology cylinders, using
our extension of the Johnson homomorphism? See Remark \ref{rem.fil}.
\end{question}

A positive answer to the above question would imply that the full tree-level
part of the theory of \fti s on 3-manifolds is given by our extension
of the Johnson homomorphism to homology cylinders.

\begin{question}
If $(\sti )-M\in\F^Y_n (\sti )$, is $M\in\H_{g,1}(n)\otimes\BQ$? See Remark
\ref{rem.fil}. When $g=0$ the answer is yes---see \cite{GL2}.  Are the
$\mu$-invariants of homology cylinders (see Remark \ref{rem.mu})
finite-type in the Goussarov-Habiro sense?
\end{question}

We now consider the group $\H^c_{g,1}$ of homology cobordism classes of
homology cylinders defined in Remark \ref{rem.hc}. The subgroup
$\H^c_{g,1}[2]$ (see Proposition \ref{prop.DA}) is the analogue of the
Torelli group, which we denote $\T\H^c_{g,1}$. We can consider the lower
central series filtration $(\T\H^c_{g,1})_n$ and, just as in the case of
the Torelli group (see \cite{Mo2}), we have
$(\T\H^c_{g,1})_n\sub\H^c_{g,1}[n+1]$. Recall also that Hain \cite{Hain}
proved that, over $\BQ$, the associated graded Lie algebra of the lower
central
series filtration of the Torelli group maps onto the associated graded Lie
algebra of the weight filtration
of the Torelli group, and that these two filtrations are known to differ
in degree $2$ by a factor of $\BQ$, \cite{Hain,Mo3}. Whether they differ
in degrees other than $2$ is an interesting question.

\begin{question}
What is the relation between the filtrations $ (\T\H^c_{g,1})_n$ and $
H^c_{g,1}[n]$?
\end{question}

Notice that the answer to above question is not known for the group
of concordance classes of string-links, see \cite{HM}, but it is
known (in the positive) for the group of homotopy classes of string-links,
see \cite{HL}, and for the pure braid group, see \cite{Kh}.

\begin{question}
We have $(\T_{g,1})_n\sub\H^c_{g,1}(n)$ (see Remark \ref{rem.fil}) and
obviously $(\T_{g,1})_n\sub (\T\H^c_{g,1})_n$. What is the relation between
the filtrations $ \H^c_{g,1}(n)$ and $ (\T\H^c_{g,1})_n$?
\end{question}

We now consider the center $\mathcal Z (\H^c_{g,1})$ of $\H^c_{g,1}$. This
contains, at least, the group $\theta^H$ of homology cobordism classes of
homology $3$-spheres (see Remark \ref{rem.hc}). Furthermore $\mathcal Z
(\H^c_{g,1})$ contains also the element $\tau_{\bd}$ defined by a Dehn
twist about the boundary of $\S_{g,1}$. For the mapping class group
$\Gamma_g$ it seems to be true (according to J. Birman and C. McMullen)
that the center is trivial, at least for large $g$.

\begin{question} Determine $\mathcal Z (\H^c_{g,1})$. Is it generated by
$\theta^H$ and $\tau_{\bd}$?
\end{question}

We now consider the subgroup $ \H^c_{g,1}[\infty]\sub \H^c_{g,1}$. This
contains $\theta^H\times \mathcal S_g^c [\infty ]$, where $ \mathcal S_g^c
[\infty ]$ denotes the subgroup of the string-link concordance group $
\mathcal S_g^c$ consisting of string links with vanishing $\mu$-invariants
(see Remarks \ref{rem.hc} and \ref{rem.sl}).  $\mathcal S_g^c [\infty ]$
contains, for example, all boundary string links and, in particular, the
knot concordance group, which is an abelian group of infinite rank
(see\cite{Le1}).

\begin{question}\lbl{q.inf}
Determine  $ \H^c_{g,1}[\infty]$. Is it equal to $\theta^H\times\mathcal
S_g^c [\infty ]$?
\end{question}

\begin{question}
Is  $ \mathcal Z (\H^c_{g,1})/
\mathcal Z (\H^c_{g,1})\cap \H^c_{g,1}[\infty]$ the infinite cyclic group
generated by $\tau_\partial$?
\end{question}

In \cite{Mo3}, Morita calculated the
symplectic invariant part $\LLc(H)^{\mathfrak s\mathfrak p}$ of $\LLc(H)$
in terms of a beautiful space of chord diagrams.
The group $ \mathcal Z (\H^c_{g,1})/(
\mathcal Z (\H^c_{g,1})\cap \H^c_{g,1}[\infty])$ is closely related
to $\LLc(H)^{\mathfrak s\mathfrak p}$, since the image under the map
$\sigma_n$ of Theorem \ref{th.alf} of an element in $
\mathcal Z (\H^c_{g,1})\cap \H^c_{g,1}[n]$ lies in
$\LLc_n(H)^{\mathfrak s\mathfrak p}$.
Thus any element of $ \mathcal Z (\H^c_{g,1})/
\mathcal Z (\H^c_{g,1})\cap \H^c_{g,1}[\infty]$ provides a geometric
construction of an element of $\LLc(H)^{\mathfrak s\mathfrak p}$.

The following question is important to the philosophical notion of
finite type.

\begin{question}
Is  $\H^c_{g,1}$ finitely-generated?
Is its abelianization finitely-generated?
Note  that both the mapping
class group $\Ga_{g,1}$ and the Torelli group $\T_{g,1}$ are
finitely-generated.
Note also that $\mathcal S_g^c$ and $ \theta^H$ are infinitely-generated
abelian (see
\cite{F}) and, as for the analogous question for string-links, since the knot
concordance group has infinite rank, the
abelianization of the string-link concordance group (on any number of strings)
has infinite rank.
\end{question}

\begin{question} Let $\H_g$ denote the semigroup of homology cylinders over
the {\em closed } surface $\Sigma_g$ of genus $g$. The kernel of the
obvious epimorphism $\H_{g,1}\to\H_g$ is related to concordance classes of
framed proper arcs in $I\times \Sigma_g$. Describe this more explicitly
and consider also the homology cobordism groups $\H_g^c$.
\end{question}

\subsection{Note}
The present paper was completed in 1999, and its follow-up by the second 
author appeared in \cite{Le3}. 


\ifx\undefined\bysame
	\newcommand{\bysame}{\leavevmode\hbox
to3em{\hrulefill}\,}
\fi

\end{document}